\documentclass[10pt,a4paper,twoside]{article}
\usepackage{amsmath} \usepackage{amssymb}
\usepackage{latexsym} \usepackage{mathrsfs} 
\usepackage{fancyhdr} 
\usepackage{mathrsfs}
\usepackage{mathptmx}

\newtheorem{definition}{\bf Definition}[section]

\newtheorem{lemma}[definition]{\bf Lemma}
\newtheorem{theorem}[definition]{\bf Theorem}
\newtheorem{proposition}[definition]{\bf Proposition}

\newcommand{\qed}{\rightline{$\Box $}}

\newcommand{\fil}{{\rm Fil}}
\newcommand{\nabs}{\nabla_{\rm abs}} 

\newcommand{\setmin}{\! \smallsetminus \!}
\newcommand{\Hol}{{\cal H}} 
\newcommand{\R}{{\mathbb R}} \newcommand{\C}{{\mathbb C}} 
\newcommand{\N}{{\mathbb N}}
\newcommand{\Z}{{\mathbb Z}} 
 \newcommand{\cc}{{\cal C}}

\newcommand{\MM}{{\cal M}}
\renewcommand{\S}{{\mathfrak S}} 
\renewcommand{\o}{{\scriptstyle{\cal O}}}
\newcommand{\sD}{\ast D}
\newcommand{\vi}{\varphi} \newcommand{\ve}{\varepsilon} 
\newcommand{\bve}{{\cal L}}
\newcommand{\vt}{\vartheta} 
\newcommand{\vve}{\varepsilon^\vee}
\newcommand{\Int}{\int\limits}
\newcommand{\ot}{\leftarrow} \newcommand{\lot}{\longleftarrow}
\newcommand{\lto}{\longrightarrow}
\newcommand{\Crd}{{\cal C}^{rd}} \newcommand{\Krd}{{\cal K}^{rd}} 
\newcommand{\Hrd}{H^{rd}}
\newcommand{\Hdr}{H_{dR}} 

\newcommand{\jE}{j_\ast (L|_U)} \newcommand{\Eess}{L_{\rm ess}}
\newcommand{\Emero}{L_{\rm mero}}
\newcommand{\Enu}{L_\nu} 
\newcommand{\deR}[1]{{\rm DR}(#1)}
\renewcommand{\mod}[2]{\raisebox{.4ex}{$#1$} \!\! \bigm/ \!\!
  \raisebox{-.4ex}{$#2$}}  
\renewcommand{\k}{\kappa} \renewcommand{\l}{\lambda}
\newcommand{\M}{{\mathfrak M}}
\renewcommand{\P}{{\mathbb P}}
 
\renewcommand{\O}{{\cal O}} 
\newcommand{\pf}{{\bf Proof: }}

\makeatletter
\renewcommand{\section}{\@startsection{section}{1}{\z@}%
{-\baselineskip}{1\baselineskip}{\large\bf}}
\renewcommand{\subsection}{\@startsection{subsection}{2}{\z@}%
{-\baselineskip}{0.5\baselineskip}{\normalsize\bf}}
\renewcommand{\subsubsection}{\@startsection{subsubsection}{3}{\z@}%
{-\baselineskip}{0.5\baselineskip}{\normalsize\bf}}
\makeatother

\numberwithin{equation}{section} 

\begin{document}

\pagestyle{plain}
\headheight12pt \renewcommand{\headrulewidth}{0.1pt} \fancyhf{}
\fancyhead[ER]{Periods for rank 1 irregular singular
  connections on surfaces}
\fancyhead[EL]{\thepage} \fancyhead[OL]{Periods for rank 1 irregular singular
  connections on surfaces} \fancyhead[OR]{\thepage}

\centerline{\Large\bf Periods for rank 1 irregular singular
  connections}
\centerline{\Large\bf on surfaces}

\bigskip

\centerline{{\bf  Marco Hien}\footnote{NWF I -- Mathematik, Universit{\"a}t
  Regensburg, 93040 Regensburg, Germany, email:
  marco.hien@mathematik.uni-regensburg.de}}

\bigskip\noindent
\centerline{\parbox[b]{14cm}{{\bf Abstract:}
  We define a period pairing for any flat, irregular singular, rank one
  connection satisfying a technical condition
  regarding its stationary set on a complex surface between de Rham
  cohomology of the connection and a modified singular homology, the
  rapid decay homology. We prove that this gives a perfect duality.}}

\section{Introduction}

Let $X $ be a smooth quasi-projective algebraic variety over the
complex numbers and $E $ be a vector bundle on $X $ equipped with an
integrable connection
$$
\nabla: E \to E \otimes_{\O_X} \Omega^1_{X|\C} \ .
$$
Its de Rham cohomology $\Hdr^\ast(X; E, \nabla) $ is defined as the
hypercohomology of the complex
$$
0 \lto E \stackrel{\nabla}{\lto} E \otimes_{\O_X} \Omega^1_{X|\C}
\stackrel{\nabla}{\lto} E \otimes_{\O_X} \Omega^2_{X|\C}
\stackrel{\nabla}{\lto} \ldots \stackrel{\nabla}{\lto} E
\otimes_{\O_X} \Omega^{\dim X}_{X|\C} \lto 0 \ ,
$$
where $\Omega^p_{X|\C} $ denotes the sheaf of K\"ahler $p $-forms on
$X $.

In addition to these data, we can consider the analytic manifold
$X_{an} $ associated to $X $ as well as the associated analytic vector
bundle $E_{an} $ with holomorphic connection $\nabla_{an} $. The
hypercohomology of the resulting complex of sheaves of holomorphic
forms with values in $E $ gives the analytic de Rham cohomology. If $X
$ is projective, it follows from the Poincar\'e Lemma and Serre's GAGA,
that the algebraic and the analytic de Rham cohomology coincide.
Equivalently, if ${\cal E}^\vee $ denotes the local system of
solutions of the dual connection $\nabla_{an}^\vee $ on the dual
bundle $E_{an}^\vee $, integration defines a perfect pairing
\begin{equation} \label{eq:dual} \Hdr^\ast(X; E, \nabla) \times
  H_\ast(X_{an}, {\cal E}^\vee) \to \C
\end{equation}
between algebraic de Rham cohomology and singular homology with values
in the local system ${\cal E}^\vee $.

If we start with a quasi-projective variety $U $, which we consider to
be compactified by a projective variety $X $, the situation
is more complicated. Let $D:= X \smallsetminus U $ denote
the complement which we 
assume to be a normal crossing divisor. In \cite{deligne}, P. Deligne
introduces the condition for a connection to be regular singular along $D
$ generalizing the well-known property for linear
differential operators in one variable (Fuchs condition) and proves the comparison
isomorphism and hence the perfect duality of \eqref{eq:dual} under
this assumption (\cite{deligne}, Th\'eor\`eme II.6.2).

In the irregular singular case, the period pairing as in
\eqref{eq:dual} is no longer perfect. The appropriate generalization
is known in dimension one only (cp. \cite{b-e}).  On curves, S. Bloch and H. Esnault
define a modified homology, the {\em rapid decay homology groups}
$\Hrd_\ast(X_{an}; E_{an}, \nabla_{an}) $ and obtain a perfect duality
$$
\Hdr^\ast(U; E, \nabla) \times \Hrd_\ast(X_{an}, E_{an}^\vee,
\nabla_{an}^\vee) \to \C \ ,
$$
given by integration. The resulting periods are interesting objects by
themselves (the integral representations of the classical
Bessel-functions, Gamma-function and confluent hypergeometric
functions arise in this way as periods of irregular singular rank one
connections on curves) and are mysteriously related to ramification
data for certain wildly ramified $\ell $-adic sheaves on curves over a
finite field (see e.g. \cite{terasoma}).

In the present paper, we want to start the investigation of the
higher-dimensio\-nal case by studying the period pairing for irregular
rank one connections on complex surfaces. We work entirely in the
analytic topology, the algebraic aspect we originally have in mind
will be mirrored by looking at an integrable connection on the smooth
analytic manifold $X $ which is meromorphic along the normal crossing
divisor $D $ at infinity, i .e. the connection is given as
$$
\nabla: E(\ast D) \to E \otimes_{\O_X} \Omega^1_{X|\C}(\ast D) \ ,
$$
where we use the usual notation $E(\ast D) $ for the sheaf of sections
of $E $ meromorphic along $D $ and we skip the subscript {\em an} in
the following.

Furthermore, we restrict ourselves to the case of line bundles $L $ with
irregular singular connections, which we assume to be {\em good with
  respect to the divisor} $D $, which is defined as follows: Consider the
formal connection $\widehat{L} := L \otimes \O_{\widehat{X|Y}}(\sD) $,
where $\O_{\widehat{X|Y}} $ denotes the formal completion of $\O_X $
with respect to the stratum $Y $ of $D $ considered, namely $Y $ being
a smooth component of $D $ or a crossing-point. By a standard
argument, $\widehat{L} $ is locally isomorphic to the formal completion of a
connection of the form $e^\alpha \otimes R $, with $\alpha \in
\O_X(\ast D) $ (cp. \cite{sabbah1}, Proposition III.2.2.1), where
$e^\alpha $ denotes the connection on the
trivial bundle $\O_X $ given by $\nabla 1 = d \alpha $ (such that the
local solutions are of the form $e^\alpha $), and $R $ is a regular
singular connection.
\begin{definition} \label{def:good}
The connection $\nabla $ is {\bf good with respect to} $D $, if its
formal completion locally is isomorphic to the formal completion of
$e^\alpha \otimes R $ with a regular singular $R $ and a local section
$\alpha \in \O_X(\ast D) $, such that the divisor $(\alpha) $ of
$\alpha $ is contained in $D $ and negative.
\end{definition}
We thus exclude examples like $\alpha= x_1^{-m_1}-x_2^{-m_2}
$. Remark, that the definition means that a good connection $\nabla $
has a local formal presentation as above, locally at the point $x=0 $
with coordinates
$x_1, x_2 $ with $D=\{ x_1 x_2 =0 \} $ or $D=\{ x_1=0 \} $, such that
$$
\alpha=x_1^{-m_1} x_2^{-m_2} \cdot u(x)
$$
with $u(0) \neq 0 $. We will always assume this. Note that any rank
one connection will become good after a finite number of point
blow-ups centered at points on $D $. We will come back to the more
general situation (for higher rank connections) in a follwing paper.

We generalize the notion of rapid decay homology groups and prove that
the resulting period pairing between the meromorphic de Rham
cohomology and rapid decay homology is perfect:
\begin{theorem}
  Let $L $ be a line bundle on a smooth projective complex surface $X $. If
  $\nabla $ is an integrable connection which is meromorphic 
  along the normal crossing divisor $D \subset X $ and good with
  respect to $D $, the period pairing
$$
\Hdr^\ast(X \smallsetminus D; L, \nabla) \times \Hrd_\ast(X; L^\vee,
\nabla^\vee) \lto \C
$$
is a perfect duality.
\end{theorem}
\pagestyle{fancy}

In dimension one, the Levelt-Turrittin
theorem and the theory of Stokes structures allows to reduce the
higher rank case to the case
of irregular singular line bundles (cp. \cite{b-e}). On surfaces,
there are analogous partial results for higher rank connections due to
C. Sabbah, e.g. an analogue to the Levelt-Turrittin theorem in the
case of rank less than or equal to $5 $. However, there are subtle
differences between the one- and the two-dimensional situation, mainly
concerning the non-good situation, with
interesting consequences regarding the period pairing as well. We will
come back to this in a subsequent paper.

Additionally, it turns out to be very difficult to give explicit
examples of flat meromorphic connections of higher rank due to the
integrability condition imposed. Locally, a rank $r $ connection is
given by its connection matrix $A = A_1 dx_1 + A_2 dx_2 $ with $r
\times r $-matrices $A_i $ having meromorphic functions as entries. The
integrability condition reads as
$$
\frac{\partial A_2}{\partial x_1} - \frac{\partial A_1}{\partial x_2}
= [ A_1, A_2 ] \ ,
$$
and it is difficult to find explicit (non-trivial) meromorphic
solutions to this equation (cp. the corresponding remark in
\cite{sabbah1}, p. 2). It is however possible to construct
higher-rank examples by functoriality (in the category of ${\cal D}
$-modules), e.g. by pushing forward an irregular singular rank one
connection $(L, \nabla) $ on $X $ along some map $f:X \to Y $. The
resulting Gau\ss-Manin connection lives on the higher direct image
$R^pf_{\ast,dR} (L, \nabla) $ of the de Rham complex on $X $. We will
give an example, a two-dimensional generalization of the confluent
hypergeometric connection, at the end of this introduction,
leading to an explicit meromorphic rank $3 $ connection on $\P^2 $.

Examples involving line bundles, however, occur in a natural way in
the framework of special functions, more precisely the higher
dimensional generalizations of well-studied special functions, such as
the generalized hypergeometrics in the sense of Gelfand and Aomoto and
their confluent variants (cp. \cite{haraoka}, \cite{kihata}). For
example, any closed meromorphic $1 $-form $\omega $ on $X $ with
poles along the normal crossing divisor $D $ gives rise to a singular
connection on the trivial line bundle given by
$$
\nabla = d + \omega \wedge \ .
$$
If the $1 $-form $\omega $ depends on additional parameters, say $z
\in Z $, and if $\omega $ is closed as a $1
$-form on the product $X \times Z $ (and hence induces an integrable
connection on the trivial line bundle $\O_{X \times Z} $), the
resulting periods for the relative connection on $X $ satisfy the
Gau{\ss}-Manin connection on $Z $. Bessel-functions and more generally
confluent hypergeometric functions occur in this way. We want to
illustrate this construction with an example.

\subsection{An example: two-dimensional confluent hypergeometrics}

We fix $a,b,c,\alpha \in \C $, with $a,b,c \not\in \Z $ satisfying
$a+b+c=-3 $, as well as additional parameters $x,y $. Consider the
connection on the trivial line bundle $\O_{\P^2} $ on $\P^2 $, which
in affine coordinates $[1,u_1,u_2] \in \P^2 $ reads as
$$
\nabla = d + \big(\frac{a}{u_1} + \frac{c}{1+u_1+u_2}+ \alpha x
\big) du_1 + \big( \frac{b}{u_2}+ \frac{c}{1+u_1+u_2} + \alpha y
\big) du_2 \ .
$$
Solutions are given by the various branches of
$$
U(u_1,u_2,x,y) := u_1^a \cdot u_2^b \cdot (1+u_1+u_2)^c \cdot \exp(\alpha(x
u_1+y u_2)) \ .
$$
The connection above can be written as $\nabla=d+ {\rm dlog}_u U $,
where ${\rm dlog}_u $ denotes the logarithmic derivation with respect
to the coordinates $u=(u_1, u_2) $, i.e. ${\rm dlog}_u U= {\rm
  dlog}_{u_1} U \, du_1 + {\rm dlog}_{u_2} U \, du_2 $.

If we interpret the parameters $x, y $ as coordinates in the affine
space and compactify with a projective plane at infinity, i.e. we
read $(x, y) $ as the point $[1,x,y] \in \P^2 $, we obtain an
integrable meromorphic connection $\nabla $ on the trivial line bundle
over $\P^2 \times \P^2 $, namely in affine coordinates
$$
\nabs = d + {\rm dlog}_{(u,x,y)} U = d + {\rm dlog}_u U + \alpha u_1
dx + \alpha u_2 dy \ .
$$
Let $D := \{ u_1u_2=0 \} \cup \{ 1+u_1+u_2=0 \} \cup \{ [0,t_1,t_2]
\in \P^2 \} $ and $X := \P^2 \smallsetminus D $. On one-forms, the
connection $\nabla: \Omega^1_X(\ast D) \to \Omega^2_X(\ast D) $ reads
as $\nabla(f du_1 + g du_2) = $
$$
=\big( \frac{\partial g}{\partial u_1} -
\frac{\partial f}{\partial u_2} + g \cdot (\frac{a}{u_1} +
\frac{c}{1+u_1+u_2} + \alpha x) - f \cdot ( \frac{b}{u_2} +
\frac{c}{1+u_1+u_2} + \alpha y) \big) du \ ,
$$
where we abbreviate $du := du_1 \wedge du_2 $. We claim, that $\dim
\Hdr^2(X, \O_{X}(\ast D), \nabla) = 3 $, a basis is given by the de
Rham classes of $du $, $u_1 du $, $u_2 du $. In order to understand
this, we consider the following equalities in $\Hdr^2 := \Hdr^2(X;
\O_{X}(\ast D), \nabla) $. First, we have $0 \equiv \nabla \big( u_2
(1+u_1+u_2) du_1 \big) = $
\begin{equation} \label{eq:GMI} 
=\big( \alpha x u_1^2 + \alpha x u_1u_2 + (\alpha x - (1+b))
  u_1 + (1+a) u_2 + (1+a) \big) du \ ,
\end{equation}
as well as $0 \equiv \nabla \big( u_1 (1+u_1+u_2) du_2
  \big) = $
\begin{equation} \label{eq:GMII}  
= - \big( \alpha y u_2^2 + \alpha y u_1u_2 + (1+b) u_1 +
  (\alpha y - (1+a)) u_2 + (1+b) \big) du \ .
\end{equation}
Additionally, calculating
\begin{equation} \label{eq:nab1} \nabla \big( (\alpha xu_1 +
  (1+a))u_2(1+u_1+u_2)du_1 + (\alpha yu_2 + (1+b)) u_1(1+u_1+u_2)du_2)
\end{equation}
induces the following equivalence in $\Hdr^2 $:
\begin{equation} \label{eq:GMIII} \alpha u_1u_2 du = -
  \frac{1+b}{y-x} u_1 du + \frac{1+a}{y-x} u_2 du \ .
\end{equation}
Evaluating higher powers $u_1^ku_2^l(1+u_1+u_2)du_i $, $i=1,2 $, and
similar linear combinations as in \eqref{eq:nab1}, easily proves the
claim about $\Hdr^2 $.

\noindent
The period pairing
$$
\Hdr^2(X; \O_{X}(\ast D), \nabla) \times \Hrd_2(\P^2; \O_{X}(\ast D),
\nabla) \to \C \ ,
$$
which we are going to define in the next section, produces the
following kind of generalized confluent hypergeometric functions on $X
$:
$$
F(x,y) = \Int_c u_1^a \cdot u_2^b \cdot (1+u_1+u_2)^c \cdot
\exp(\alpha (x u_1+y u_2)) \omega
$$
with $\omega $ being either $du $, $u_1 du $ or $u_2 du $. How the
topological $2 $-chain $c $ in $\P^2 $ has to be chosen, will be the
main point in the definition of the rapid decay homology $\Hrd_\ast $.

These confluent hypergeometric functions $F $ again satisfy another
partial differential equation, namely the Gau\ss-Manin equation on
$\Hdr^2 $ derived from $\nabla $. The latter is defined as follows.
Let $Z=\P^2 $ denote the space for the parameters $(x,y) $ and let
$\pi: X \times Z \to Z $ be the projection. We will keep the affine
coordinates $u=(u_1,u_2) $ for the points in $X $ and $(x,y) $ in $Z
$. There is a filtration on $\Omega^\ast_{X \times Z |\C} $ given by
$$
\fil^i \Omega^\ast_{X \times Z|\C} := {\rm im}(\pi^\ast
\Omega^i_{Z|\C} \otimes \Omega^{\ast-i}_{X \times Z|\C}) \ .
$$
Now, the associated graded object fulfills ${\rm gr}^i \cong \pi^\ast
\Omega^i_{Z|\C} \otimes \Omega^{\ast-i}_{X \times Z|Z} $, especially
${\rm gr}^0 \cong \Omega^i_{X \times Z|Z} $, the sheaf of relative
differential forms, on which the original connection $\nabla $
canonically lives. The short exact sequence $0 \to {\rm gr}^1 \to
\fil^0/\fil^2 \to {\rm gr}^0 \to 0 $ induces the short exact sequence
of de Rham complexes
\begin{equation} \label{eq:GMdiag} 0 \to (\Omega^{\ast-1}_{X \times
    Z|Z} \otimes \pi^\ast \Omega^1_{Z|\C}, \nabs \otimes 1) \to
  \big( \mod{\Omega^\ast_{X \times Z|\C}}{\fil^2}, \nabs \big) \to
  (\Omega^\ast_{X \times Z|Z}, \nabla) \to 0 \ .
\end{equation}
The Gau\ss-Manin connection is by definition the connecting morphism
of the associated long exact sequence of the higher direct images, in
our case $\nabla_{\rm GM}: \Hdr^2 \to \Hdr^2 \otimes \Omega^1_{Z|\C}
$ (note that $\Hdr^2 = R^2\pi_\ast(\Omega^\ast_{X
  \times Z|Z}, \nabla) $). Chasing the diagram \eqref{eq:GMdiag} gives
$$
\nabla_{\rm GM} (u_1^k u_2^ldu) = \alpha dx \otimes u_1^{k+1}u_2^l du
+ \alpha dy \otimes u_1^k u_2^{l+1} du \ .
$$
Applying \eqref{eq:GMI}, \eqref{eq:GMII} and \eqref{eq:GMIII}, we
obtain the connection matrix $\Phi_{\rm GM} $ with respect to the
basis $du $, $u_1 du $, $u_2 du $ of $\Hdr^2 $:
\begin{equation} \label{eq:GM} \Phi_{\rm GM} = \left(
    \begin{array}{c|c|c} 0 & -\frac{1+a}{x} & 0 \\[0.2cm]
      \alpha & \frac{(1+b)y}{x(y-x)} - \alpha & -\frac{1+b}{y-x} \\[0.2cm]
      0 & -\frac{(1+a)y}{x(y-x)} & \frac{1+a}{y-x} \end{array} \right)
  \, dx + \left(
    \begin{array}{c|c|c} 0 & 0 & -\frac{1+b}{y} \\[0.2cm]
      0 & - \frac{1+b}{y-x} & \frac{(1+b)x}{y(y-x)} \\[0.2cm]
      \alpha & \frac{1+a}{y-x} & - \frac{(1+a)x}{y(y-x)} - \alpha
    \end{array} \right) \, dy \ ,
\end{equation}
being one of the rather rare explicit examples of an integrable higher
rank connection on a surface.

\medskip We remark, that there is a theory of generalized confluent
hypergeometrics on the space $Z_{r+1,n+1} $ of complex $(r+1) \times
(n+1) $-matrices of full rank for given $1 \le r <n $, defined in
\cite{kihata}. The starting point again is a connection
of the form $\nabla= d+ {\rm dlog} \, U $ for a certain class of
multi-valued functions $U: \P^r \times Z_{r+1, n+1} \to \C $ with
well-defined logarithmic derivative. Actually, one also fixes a
composition $\lambda=(1+\lambda_0, \ldots, 1+ \lambda_l) $ of $n+1 $,
i.e. $\sum (1+\lambda_i)=n+1 $ and requires for the solutions to be
invariant under the left and right action $ {\rm GL}_{r+1}(\C) \times
Z_{r+1,n+1} \times H_\lambda \to Z_{r+1,n+1} $, where $H_\lambda
\subset {\rm GL}_{n+1}(\C) $ denotes the subgroups of all block
diagonal matrices with $l+1 $ blocks consisting of upper triagonal
matrices with the constant entry $h_i^{(k)} $ along the $i $-th
upper diagonal for $i=0, \ldots, \lambda_k $ (the entry on the main
diagonal $h_0^{(k)} \neq 0$).
Our example above corresponds to the choices $r=2 $, $n=4 $ and
$\lambda=(2,1,1,1) $ and the restriction to the subspace of all
matrices
$$
\left( \begin{array}{cc|c|c|c} 1& 0 & 0 & 0 & 1 \\ 0 & x & 1 & 0 & 1
    \\ 0 & y & 0 & 1 & 1 \end{array} \right) \in Z_{3,5} \ ,
$$
which parameterize the {\em generic stratum} of the double quotient
${\rm GL}_3(\C) \backslash Z_{3,5} / H_\lambda $ (we refer to
\cite{kihata} and \cite{haraoka} for further details). 

\medskip\noindent The paper is organized as follows. In section 2, we
define the rapid decay homology groups for line bundles on $X $ and
its pairing with the meromorphic de Rham cohomology. Afterwards, we
reduce the problem of perfectness of the pairing to local questions
according to the canonical stratification of the normal crossing
divisor $D $ into crossing points and smooth components. Perfectness
of the resulting local pairings is proved in several steps in section
4, completing the proof of the main theorem.

\section{Rapid decay homology and the pairing with de Rham cohomology}

\subsection{Rapid decay homology}

Let $X $ be a $n $-dimensional smooth projective complex manifold and $D \subset
X $ a divisor with normal crossings (i.e. in suitable local
coordinates $z_1, \ldots, z_n $ it is of the form $D=\{ z_1 \cdots
z_k=0 \} $ for some $1 \leq k \leq n $ --- such coordinates will be
called {\em good} w.r.t $D $). We further consider a line bundle $L $
over $X $ together with an integrable meromorphic connection on
$X\setmin D $ with possibly irregular singularities at $D $. In the
usual notation $L(\sD) $ for the sheaf of local sections in $L $
meromorphic along $D $, the connection then reads as $\nabla: L(\sD)
\lto L \otimes \Omega^1(\sD) $. The rank one local system of
horizontal sections in $L $ on the complement $U:=X \setmin D $ will
be denoted by
$$
\bve := \ker(L|_U \stackrel{\nabla}{\to} L|_U \otimes \Omega_U^1)
\subset L|_U \ .
$$
The dual connection $\nabla^\vee $ on $L^\vee $ is characterized by $d
<e,\vi> = < \nabla e, \vi > + < e, \nabla^\vee \vi > $ for local
sections $e $ of $L $ and $\vi $ of $L^\vee $. Let $\bve^\vee $ denote
the corresponding local system.

We assume that the complement $U:= X\setmin D $ is Stein, which we can
always obtain by joining additional hypersurfaces, where the
connection is not singular at all, to $D $. These do not affect our
procedure. Then the de Rham cohomology is given by the cohomology of
the global sections in $U $, i.e.
$$
H^p_{dR}(U; L(\ast D),\nabla) = H^p( \ldots \to \Gamma(L \otimes
\Omega^q(\sD)) \stackrel{\nabla}{\to} \Gamma(L \otimes
\Omega^{q+1}(\sD)) \to \ldots) \ .
$$

The homology we are going to define will be a generalization of the
usual notion of singular homology with coefficients in a local system
${\cal V} $, where one considers chain complexes built from pairs of a
topological chain together with an element of the stalk of ${\cal V} $
at the barycentre of the chain.  In our situation, we will have to
allow the topological chains to be able to meet the divisor $D $,
where the local system $\bve := (L|_U)^\nabla $ is not defined. To
this end, we will make the following definition:
\begin{definition}
  For any $x \in X $, we define the {\bf stalk} $\bve_x $ of $\bve $
  to be the usual stalk of $\bve $ if $x \in U $ and to be the
  coinvariants
  $$
  \bve_x := (\bve_y)_{\pi_1(U,y)}:= \mod{\bve_y}{\{ v- \sigma v \mid v
    \in \bve_y, \sigma \in \pi_1(U,y) \}} \mbox{\qquad for } x \in D \
  ,
  $$
  where $y \in U $ is any point near $x \in D $ and we are taking
  coinvariants w.r.t. the local monodromy action $\pi_1(U,y) \to {\rm
    Aut}(\bve_y) $.
\end{definition}

Now, we can define the notion of rapidly decaying topological chains
in analogy to the definition by S. Bloch and H. Esnault in \cite{b-e}
on curves. In the following, we denote by $\Delta^p $ the standard $p
$-simplex with barycentre $b \in \Delta^p $ and we call a function
$f:Y \to \C $ from any subset $Y \subset U $ of the open complex
manifold $U $ analytic, if it is the restriction of an analytic
function on an open neighborhood of $Y $.

\begin{definition}
  A {\bf rapid decay $p $-chain} is a pair $(c, \ve) $ consisting of a
  continuous map $c: \Delta^p \to X $, such that the pre-image
  $c^{-1}(D) $ is a union of complete subsimplices of $\Delta^p $, and
  an element $\ve \in \bve_{c(b)} $. If $c(\Delta^p) \not\subset D $,
  we require that $\ve \in \bve_{c(b)} \subset L_{c(b)} $ is {\rm
    rapidly decaying} in the following sense:
  
  For any $y \in c(\Delta^p) \cap D $, let $e $ denote a local
  trivialization of $L(\ast D) $ and $z_1, \ldots, z_n $ local
  coordinates at $y $ such that locally $D = \{ z_1 \cdots z_k = 0 \}
  $ and that the coordinates of $y $ fulfill $y_1 = \ldots = y_k = 0
  $. With respect to the trivialization $e $ of $L $ restricted to $U $,
  $\ve $ becomes an analytic function
  $$
  f:= e^\ast \ve: c(\Delta^p)\setmin D \to \C \ , \ (z_1, \ldots, z_n)
  \mapsto f(z) \ .
  $$
  We require that this function has rapid decay at $y $, i.e.  that
  for all $N \in\N^k $ there is a $C_N > 0 $ such that
  $$
  |f(z)| \le C_N \cdot |z_1|^{N_1} \cdots |z_k|^{N_k}
  $$
  for all $z \in c(\Delta^p) \setmin D $ with small $|z_1|, \ldots,
  |z_k| $.
\end{definition}

We stress that we do not impose any decay condition on pairs $(c, \ve)
$ with $c(\Delta^p) \subset D $; nevertheless we call those pairs
rapidly decaying as well.

Now, let $S^D_p X $ be the free vector-space over all singular chains
$\Delta^p \to X $ meeting $D $ only in full subsimplices and let
$\Krd_p(X;L) $ be the $\C $-vector space of all maps
$$
\psi: S^D_p X \to \bigsqcup_{x \in X} \bve_x ,
$$
such that $\psi(c)=0 $ for all but finitely many $c $ and that $(c,
\psi(c)) $ has rapid decay. We will write $c \otimes \ve \in
\Krd_p(X;L) $ for the element $\psi $ which takes the value $\ve $ at
$c $ and zero otherwise. We remark that any element of $\Krd_p(X;L) $
can be written as a finite sum $\sum^r_{i=1} c_i \otimes \ve_i $. The
notation $\otimes $ is justified by the fact, that it is linear in
each of the entries. There is a natural boundary map $\partial:
\Krd_p(X; L) \to \Krd_{p-1}(X; L) $, $c \otimes \ve \mapsto \sum_j
(-1)^j c_j \otimes \ve_j $, where the sum runs over the faces $c_j $
of $c $ and the elements $\ve_j $ in the stalk of $\bve $ at $c(b_j)
$, with $b_j $ being the barycentre of the $j^{\rm th} $ face of
$\Delta^p $, are given as follows:
\begin{enumerate}
\item if $c(b_j) \not\in D $, there is a unique homotopy class of
  paths from $c(b) $ to $c(b_j) $ (e.g. induced by the line $[b,b_j] $
  in $\Delta^p $) and $\ve_j $ is defined as the analytic continuation
  along a representative.
\item if $c(b_j) \in D $, we just define $\ve_j $ to be the element
  represented by $\ve $ in the corresponding stalk of coinvariants
  under monodromy.
\end{enumerate}
It can be easily seen, that $\partial \circ \partial = 0 $ and that
$\partial $ respects the support of the chains, so that we can define
the rapid decay homology as follows:

\begin{definition}
  For $\emptyset \subseteq Z \subseteq Y \subseteq X $, put
  $$
  \Crd_p(Y,Z; L) := \mod{\Krd_p(Y; L)}{\Krd_p(Z;L) + \Krd_p(D \cap Y;
    L)} \mbox{\quad (relative version)}
$$
and $\Crd_p(Y; L) := \Crd_p(Y, \emptyset; L) $ for the absolute version. The
{\bf rapid decay homology} is defined as the homology of the
corresponding complexes and denoted by $\Hrd_p(Y,Z; L) $ and
$\Hrd_p(Y; L) $ respectively.
\end{definition}

Note that we have moduled out the chains that are mapped completely
into $D $ as they will not play any role in the pairing with
meromorphic differential forms as described below. Nevertheless, one
has to include them a priori into the definition of rapid decay chains
in order to be able to define the boundary map $\partial $. We also
include $L $ into the notation to remind that the rd-homology does not
depend on the local system alone, but on the connection.

\subsection{The pairing and statement of the main result}

\noindent
Now, if we have a rapidly decaying chain $c\ \otimes \ve \in \Crd_p(X;
L^\vee) $ in the dual bundle (with the dual connection) and a
meromorphic $p $-form $\omega \in L \otimes \Omega^p(\sD) $, then the
integral $\int_c < \ve, \omega > $ converges because the rapid decay
of $\ve $ along $c $ annihilates the moderate growth of the
meromorphic $\omega $.  Let $c_t $ denote the topological chain one
gets by cutting off a small tubular neighborhood with radius $t $
around the boundary $\partial \Delta^p $ from the given topological
chain $c $. Then, for $c \otimes \ve \in \Crd_p(X; L^\vee) $ and $\eta
\in L \otimes \Omega^{p-1}(\sD) $ a meromorphic $p-1 $-form, we have
the {\bf 'limit Stokes formula'}
$$
\Int_c < \ve, \nabla \eta > = \lim_{t \to 0} \Int_{c_t} < \ve,
\nabla \eta > = \lim_{t \to 0} \Int_{\partial c_t} < \ve, \eta > =
\Int_{\partial c-D} < \ve, \eta > \,
$$
where in the last step we used that by the given growth/decay
conditions the integral over the faces of $\partial c_t $ 'converging'
against the faces of $\partial c $ contained in $D $ vanishes.

The limit Stokes formula easily shows in the standard way that
integrating a closed differential form over a given rd-cycle (i.e.
with vanishing boundary value) only depends on the de Rham class of
the differential form and the rd-homology class of the cycle. Thus, we
have:
\begin{proposition}
  Integration induces a well defined bilinear pairing
  \begin{equation} \label{eq:perfV}
  H^p_{dR}(U; L(\ast D),\nabla) \times \Hrd_p(X; L^\vee) \lto \C \ , \
  ([\omega], [c \otimes \ve]) \mapsto \Int_c < \omega, \ve > \ ,
  \end{equation}
  which we call the {\bf period pairing} of $(L, \nabla) $.
\end{proposition}

\noindent
Our main result is the following
\begin{theorem} \label{thm:main}
  On a complex surface $X $, the period pairing is perfect for any
  irregular singular line bundle which is good with respect to $D $.
\end{theorem}
Remark, that for ${\rm dim}(X)=1 $ (and arbitrary rank) the
perfectness was proved by S. Bloch and H. Esnault \cite{b-e} (in the
one-dimensional case every line bundle is good). The case
of arbitrary dimension $\dim(X) \ge 2 $ is not known in general. 

\subsection{The irregularity pairing}

Let $j:U:=X\setmin D \hookrightarrow X $ denote the inclusion. In
addition to the de Rham cohomology $\Hdr^\ast(U; L(\ast D)) $ of
meromorphic sections, we will also have to consider the de Rham
cohomology over $U $ allowing essential singularities as well. We will
denote the corresponding sheaves by
$$
\Emero := L(\sD) \mbox{\quad and \quad} \Eess := \jE \ .
$$
Note, that the de Rham cohomology of the connection induced on $\Eess
$ pairs classically with the singular homology $H^\ast(U; \bve^\vee) $
of $U $ with coefficients in the restricted local system $\bve^\vee|_U
$. Let $C_\ast(U;\bve^\vee) $ denote the corresponding singular chain
complex. Consider the short exact sequence of de Rham
complexes on $X $:
$$
0 \lto \deR{\Emero} \lto \deR\Eess \lto \deR{\mod{\Eess}{\Emero}} \lto
0
$$
as well as the following sequence of complexes of abelian groups
$$
0 \lto C_\ast(U; \bve^\vee) \lto \Crd_\ast(X; L^\vee) \lto \Crd_\ast(X,
U; L^\vee) \lto 0 \ ,
$$
whose exactness is obvious.

Now, for $\omega \in \Eess^p $ with $\nabla \omega \in \Emero^{p+1} $
and $c \otimes \ve \in \Crd_{p+1}(X) $ with $\partial c \in
C_p(X\setmin D) + C_p(D) $, we define
$$
< [\omega], [c \otimes \ve] > := \Int_c < \ve, \nabla \omega > -
\Int_{\partial c - D} < \ve, \omega > \ .
$$
This gives a well-defined pairing
$\Hdr^p(\mod{\Eess}{\Emero}) \times \Hrd_{p+1}(X, U) \to \C $, which
fits into the long exact sequences induced:
$$
\begin{array}{ccccccccc}
  \ldots \ot & \Hdr^{p+1}(\Emero) & \ot & \Hdr^p(\mod{\Eess}{\Emero}) &
  \ot & \Hdr^p(\Eess) & \ot 
  & \Hdr^p(\Emero) & \ot \ldots \\[0.3cm]
  \ldots \to & \Hrd_{p+1}(X; L^\vee) & \to & \Hrd_{p+1}(X, U; L^\vee) &
  \to & H_p(U; \bve^\vee) & \to & \Hrd_p(X; L^\vee) & \to \ldots 
\end{array} \ .
$$
We want to fix the given connection $(L(\sD), \nabla) $ and drop $L^\vee $
and $\bve^\vee $ in the notation of the rd-homology groups from now
on. We will refer to the pairing above as the {\bf irregularity
  pairing}. Note, that the de Rham complex of $\mod{\Eess}{\Emero} $
coincides up to a shift of degrees with the irregularity complex
introduced by Z. Mebkhout (cp. \cite{mebkhout}).

\section{Localization according to stratification of the divisor}
\label{sec:lespairg}

From now on, we concentrate on the case $\dim(X)=2 $. In this section,
we want to reduce the question of perfectness of the irregularity
pairing in two steps. The first one, being rather standard, reduces to
the local situation according to an open covering of the divisor. This
leaves us with the task to study two different local situations,
namely with $D $ being of the form $D=\{ x_1=0 \} $ or $D=\{ x_1x_2=0
\} $ in suitable local coordinates. In order to understand the
situation at a crossing point, we will further split the pairing into
one concentrating at the crossing point and one determining the
contribution of the connection along the two local components meeting
at $0 $.

If ${\mathscr U} $ denotes an open covering of $X $, we have two
cohomological spectral sequences, the first one considered being
$$
E_2^{p,q} = H^p({\mathscr U}, {\cal H}_{dR}^q(\mod{\Eess}{\Emero}))
\Longrightarrow \Hdr^{p+q}(X; \mod{\Eess}{\Emero}) \ ,
$$
where ${\cal H}_{dR}^q(\mod{\Eess}{\Emero}) $ denotes the presheaf $U'
\mapsto \Hdr^p(U', \mod{\Eess}{\Emero}) $. In the same way, one has a
homological spectral sequence involving the sheafified rapid decay
complex ${\mathscr C}_\ast^{rd} $. To define the latter, we contruct
the sheaf ${\mathscr C}^\ast_{rd} $ of rapid decay cochains as the
sheaf associated to the presheaf $U' \mapsto {\rm Hom} (\Crd_\ast(U',
U' \setmin D; D), \C) $ and let 
$$
{\mathscr C}_\ast^{rd}:= {\cal H}om_\C ( {\mathscr C}^\ast_{rd}, \C )
$$
be its dual sheaf. We consider the dual (cohomological) spectral sequence
$$
\widetilde{E}_2^{p,q} = H_p({\mathscr U}, {\mathscr C}_q^{rd})^\vee
\Longrightarrow \Hrd_{p+q}(X, X \setmin D; D)^\vee \ .
$$
Obviously, the pairing above induces a morphism between these two spectral
sequences. In order to prove perfectness of \eqref{eq:perfV}, we
therefore can do so for the local situation. i.e. a suitable open
covering. For small enough open $U' $, we thus have to consider two
cases, the one at a crossing point and the one at a smooth point of $D
$. Thus, we can assume, that $X $ is a small bi-disc $X=\Delta \cong
D^2 \times D^2 \subset \C^2 $ and $D $ reads as either the smooth
divisor $D=\{ x_1=0 \} $ or the union of two coordinate planes $D=\{
x_1x_2=0 \} $ with crossing point $(0,0) $.

\subsection{Distinguishing the contribution of the local strata}
\label{sec:locseq}

\subsubsection{The local (co-)homology groups at a crossing-point}

We consider the second case from above, i.e. $X=\Delta $ a small
bi-disc around the crossing-point $0 \in D= \{ x_1x_2=0 \} $ in
suitable local coordinates. We write $D=D_1 \cup D_2 $ for the two
local components, $D_\nu=\{ x_\nu=0 \} $. Let $j_\nu:\Delta \setmin
D_\nu \to \Delta $ denote the inclusion for $\nu=1,2 $. We will write
$\Enu^p $ for the subsheaf of $\Eess^p := \Eess \otimes \Omega^p $
defined as
$$
\Enu^p := j_{\nu \ast} j_\nu^\ast (L(\sD) \otimes \Omega^p) \ ,
$$
so that a local section $u $ of $\Enu^p $ is an analytic $p $-form
with values in $L $ defined on the complement of $D $, which is
meromorphic along $D_{\neq \nu} \smallsetminus 0 $ and arbitrary along
$D_\nu $. One might think of the function $1/z_2 \cdot \exp(1/z_1) $
as a typical example of an element in $L_1^0 $. The given connection
$\nabla $, being meromorphic, obviously carries $\Enu^p $ to
$\Enu^{p+1} $ and therefore gives the following complex of sheaves on
$\Delta $:
$$
\deR{\Enu}: \quad 0 \lto \Enu^0 \stackrel{\nabla}{\lto} \ldots \lto
\Enu^p \stackrel{\nabla}{\lto} \Enu^{p+1} \lto \ldots \ .
$$
In the same way, we write $\Emero^p := \Emero \otimes \Omega^p $. Now
$L_1 \cap L_2 = \Emero $ and we have the following short exact
sequence
$$
0 \ot \deR{\mod{\Eess}{L_1+L_2}} \ot \deR{\mod{\Eess}{\Emero}}
\ot \deR{\mod{L_1}{\Emero} \oplus \mod{L_2}{\Emero}} \ot 0
$$
of complexes of sheaves supported on $D $, where the map at the right
is given by $[\eta_1]+[\eta_2] \mapsto [\eta_1+\eta_2] $. It gives
rise to the long exact cohomology sequence
\begin{equation}
  \label{eq:les1}
  \ldots \ot \Hdr^{p-1}(\mod{\Eess}{L_1+L_2}) \ot \Hdr^{p-1}(\mod{\Eess}{\Emero})
  \ot \bigoplus_{i=1,2} \Hdr^{p-1}(\mod{L_i}{\Emero})
  \ot \ldots 
\end{equation}
which is the long exact sequence we are going to pair with the
corresponding sequence of rapid decay homology.

As for the rapid decay homology groups, an easy argument using the
subdivision morphism shows that $\Crd_\ast(\Delta \setmin D_1) +
\Crd_\ast(\Delta \setmin D_2) = \Crd_\ast(\Delta \setmin 0) $ and
$\Crd_\ast(\Delta \setmin D_1) \cap \Crd_\ast(\Delta \setmin D_2) =
\Crd_\ast(V \setmin D) $. It follows that the following sequence is
exact:
$$
0 \to \Crd_\ast(\Delta ,\Delta \setmin D;D) \to \bigoplus_{i=1,2} \Crd_\ast(\Delta
,\Delta \setmin D_i;D) \to \Crd_\ast(\Delta ,\Delta \setmin 0;D) \to 0 \ .
$$
Here the first map is given by $[\alpha] \mapsto [\alpha]+[\alpha] $
and the second one by $[\alpha]+[\beta] \mapsto [\alpha-\beta] $.  The
corresponding long exact homology sequence reads as
\begin{equation}
  \label{eq:rd}
  \ldots \to \Hrd_p(\Delta,\Delta\setmin D;D) \to
  \bigoplus_{i=1,2} \Hrd_p(\Delta,\Delta\setmin 
  D_i;D) \to \Hrd_p(\Delta,\Delta \setmin 0;D) \to \ldots 
\end{equation}
The members of the latter may be rewritten in the following form,
where we retract $\Delta \setmin D_\nu $ to the boundary $\partial
V_\nu $ of a small tubular neighborhood $V_\nu $ of $D_\nu $, i.e. if
$\Delta $ reads as the product of two disc $\Delta=D^2 \times D^2 $ in
local coordinates, $\partial V_1 $ will be $\partial V_1=S^1 \times
D^2 $.
\begin{lemma}
  There are natural isomorphisms $\Hrd_p(\Delta , \Delta \!\!\setmin
  D_\nu; D) \cong \Hrd_p(V_\nu, \!\partial V_\nu; D) $ for $\nu=1,2 $
  and $\Hrd_p(\Delta, \Delta \setmin 0;D) \cong \Hrd_p(\Delta,
  \partial \Delta;D) $.
\end{lemma}
\pf This follows easily by excision with the help of the subdivision
morphism, $ subd \simeq id $, on the complex of rapidly decaying
cycles, as well as the observation that one can chose retracts of e.g.
$\Delta \setmin D_\nu $ to the boundary $\partial V_\nu $ respecting
the rapid decay condition on the cycles $c \otimes \ve $.
\\ \qed 

\noindent
With these isomorphisms, the exact sequence \eqref{eq:rd} reads as
\begin{equation}
  \label{eq:les2}
  \ldots \to \Hrd_p(\Delta, \Delta \setmin D; D) \to \bigoplus_{i=1,2}
  \Hrd_p(V_1, \partial V_i;D) \to \Hrd_p(\Delta, \partial
  \Delta;D) \to \ldots 
\end{equation}
Our aim in this section is to define a pairing between the long exact
sequences \eqref{eq:les1} and \eqref{eq:les2}, that is a bilinear
pairing of their members, such that the obvious diagrams all commute.

\subsubsection{Definition of the pairings}

In perfect analogy to the case of the complex $\deR{\Eess} $ itself,
one defines for given
\begin{enumerate}
\item $\omega \in \Enu^p $ with $\nabla \omega \in \Emero^{p+1} $
  \quad and \quad ii) $c \otimes \ve \in \Crd(V_1) $ with $\partial c
  \in \Crd(\partial V_1) + \Crd(D) $:
\end{enumerate}
\begin{equation}
  \label{eq:j1pairg}
  < \omega, c\otimes \ve > := \Int_c < \ve, \nabla \omega > -
  \Int_{\partial c - D_1} < \ve, \omega > \ .
\end{equation}
We remark, that the first integral exists, as $\nabla \omega \in
\Emero^{p+1} $ is meromorphic along $D $ and therefore pairs with the
rapidly decaying $c \otimes \ve $. In the same way the second integral
converges, because $\omega \in L_1^p $ is meromorphic along $D_2 - 0
$ and $\ve $ decreases rapidly as $\partial c - D $ possibly
approaches $D_2 $.

Again, by using the limit-Stokes formula, one easily observes that
\eqref{eq:j1pairg} induces a well-defined pairing
$\Hdr^p(\mod{\Enu}{\Emero}) \times \Hrd_{p+1}(V_\nu, \partial V_\nu;
D) \lto \C $, which fits into the long exact sequence:
$$
\begin{array}{cccccccc}
  \ldots \ot \Hdr^{p+1}(\Emero) & \!\! \ot \!\! &  \Hdr^p(\mod{\Enu}{\Emero}) &
  \!\! \ot \!\! & \Hdr^p(\Enu) & \!\! \ot \!\! & \Hdr^p(\Emero) & \!\! \ot
  \ldots \\[0.5 em]
  \ldots \to \Hrd_{p+1}(\Delta; D) & \!\! \to \!\! & \Hrd_{p+1}(V_\nu,
  \partial V_\nu; D) & \!\! \to \!\! & \Hrd_p(\Delta \setmin D_\nu;
  D_{\neq \nu}) & \!\! \to 
  \!\! & \Hrd_p(\Delta; D) & \!\! \to \ldots
\end{array}
$$

Returning to the long exact sequences \eqref{eq:les1} and
\eqref{eq:les2}, which we want to relate to each other, we now prove
compatibility of the pairings defined so far:
\begin{lemma}
  \label{lemma:comp1} The maps
$$
\begin{array}{ccc}
\Hdr^p(\mod{\Eess}{\Emero}) & \lot & \Hdr^p(\mod{L_1}{\Emero}) \oplus
\Hdr^p(\mod{L_2}{\Emero}) \\[0.5 em]
\Hrd_{p+1}(\Delta, \Delta \setmin D; D) & \lto & \Hrd_{p+1}(V_1, \partial V_1; D_1)
\oplus \Hrd_{p+1}(V_2, \partial V_2; D_2)
\end{array}
$$
are compatible with the pairings between the groups
in the columns of the above diagram.
\end{lemma}
\pf Suppose we have $\eta_\nu \in \Enu^p $ with $\nabla \eta \in
\Emero^{p+1} $ and $c \otimes \ve \in \Crd_{p+1}(\Delta) $ with
$\partial (c \otimes \ve) \in \Crd_p(\Delta \setmin D)+\Crd_p(D) $.
Now, we can assume (by subdivision) that the topological chain $c $
decomposes as a sum $c=c_1 + c_2 $, such that $c_\nu \otimes \ve \in
\Crd_{p+1}(V_\nu, \partial V_\nu; D) $ with vanishing $\partial(c_\nu
\otimes \ve) $ (mod $\Crd_p(\partial V_\nu) + \Crd_p(D) $).
Restricting the section $\ve $ of $\bve^\vee $ to $c_\nu $ does
not affect its rapid decay properties.

The diagram above reads as $[\eta_1 + \eta_2] \leftarrow [\eta_1] +
[\eta_2] $ in the top row and $[c \otimes \ve] \mapsto [c_1 \otimes
\ve] + [c_2 \otimes \ve] $ in the bottom row and we have to prove that
$< [\eta_1 + \eta_2], [c \otimes \ve] > = < [\eta_1], [c_1 \otimes
\ve] > + < [\eta_2], [c_2 \otimes \ve] > $. Now, we consider the
decomposition $c=c_1 + c_2 $ and observe, that $c_2 \cap D_1 =
\emptyset $, the section $\ve $ is rapidly decreasing as $c_2 $
possibly approaches $D_2 $ and that $\eta_1 \in L_1^p $ is
meromorphic along $D_2-0 $, so that we can apply the 'limit Stokes
formula' and obtain
$$
\Int_c < \ve, \nabla \eta_1 > - \Int_{\partial c - D} < \ve, \nabla
\eta_1 > = \Int_{c_1} < \ve, \nabla \eta_1 > - \Int_{\partial c_1 -
  D} < \ve, \nabla \eta_1 > \ .
$$
The same argument applies to $c_2, \eta_2 $ instead of $c_1, \eta_1 $
proving the assertion.
\\ \qed 

\bigskip In the next step, we define a pairing
$\Hdr^{p-1}(\mod{\Eess}{L_1 + L_2}) \times \Hrd_{p+1}(\Delta,
\partial \Delta; D) \to \C $. To this end, let there be given
\begin{enumerate}
\item $\omega \in \Eess^{p-1} $ with $\nabla \omega =: \eta_1 + \eta_2
  \in L_1^p + L_2^p $ and
\item $c \otimes \ve \in \Crd_{p+1}(\Delta) $ with $\partial (c
  \otimes \ve) \in \Crd_p(\partial \Delta) + \Crd_p(D) $.
\end{enumerate}
We have $\partial \Delta = \partial V_1 \cup \partial V_2 $. By
subdivision we can decompose $c $ as the sum $c=c_1+c_2 $ with
$\partial c_\nu \cap \partial \Delta \subset
\partial V_\nu $. With these choices, we define
$\alpha(\eta_1,\eta_2,c_1 \otimes \ve, c_2 \otimes \ve) := $
\begin{equation} \label{eq:alphadef} := \sum_{i=1,2} (-1)^{i+1} \big( 
\Int_{c_i} < \ve, \nabla
  \eta_i > - \Int_{\partial c_i \cap
    \partial \Delta - D} \!\!\!\!\!\!\!< \ve, \eta_i > \big) +
  \Int_{\partial(\partial c_1 \cap \partial c_2)-D} <
  \ve, \omega > \ .
\end{equation}
A few words should be said about this definition:
\begin{enumerate}
\item Integrability of $\nabla $ ensures that $0 = \nabla \nabla
  \omega = \nabla \eta_1 + \nabla \eta_2 $ and therefore $\nabla
  \eta_1 = -\nabla \eta_2 \in L_1^{p-1} \cap L_2^{p-1} =
  \Emero^{p-1} $. This shows convergence of the first and third
  integral in \eqref{eq:alphadef}, as $\ve $ is rapidly decaying.
\item $\eta_\nu \in \Enu^p $ is meromorphic along $D_{\neq \nu}
  \smallsetminus 0 $ and therefore $< \ve, \eta_\nu > $ can be
  integrated over the chain $\partial c_\nu $ not meeting $D_\nu $,
  hence the second and fourth integral in \eqref{eq:alphadef}
  converge.
\item $\partial c_1 \cap \partial c_2 $ consist of the subsimplices
  that arose in the chosen decomposition of $c=c_1 + c_2 $, which we
  equip with the orientation induced from $c_1 $. These simplices are
  either fully contained in $\Delta \setmin D $ or in $D $. In the
  last integral in \eqref{eq:alphadef} we integrate over the first
  type completely contained in $\Delta \setmin D $, hence the integral
  is well-defined.
\end{enumerate}
\begin{lemma} \label{lemma:alphadef}
  Mapping $([\omega], [c \otimes \ve] )
  \mapsto \alpha(\eta_1,\eta_2, c_1 \!\otimes \ve, c_2 \!\otimes \ve) $
  gives a well-defined pairing
  $$
  \Hdr^{p-1}(\mod{\Eess}{L_1+L_2}) \times \Hrd_{p+1}(\Delta,
  \partial \Delta; D) \to \C \ .
  $$
\end{lemma}
\pf Let $\alpha := \alpha(\eta_1,\eta_2, c_1 \otimes \ve,
  c_2 \otimes
\ve) $ and let $c_1^t, c_2^t $ denote the simplices if we cut off a
small tube of radius $t $ from $D $, so that $c_i^t $ converges to
$c_i $ for $t \to 0 $. Then $\alpha $ is the limit $t \to 0 $ of the
sum of the integrals as in \eqref{eq:alphadef} the cut-off $c_i^t $
instead of $c_i $. We decompose $\partial c_i^t = \gamma_i^t + \nu_i^t
+ \zeta^t $, where $\gamma_i^t \subset \partial \Delta $, $\zeta^t $
denotes the simplices of $\partial c_1^t \cap \partial c_2^t $ not
running completely into $D $ and $\nu_i^t $ the simplices of $\partial
c_i^t $ running completely into $D $ for $t \to 0 $, all of them taken
with their natural orientations induced from $c $ (especially $\zeta^t
$ carries different orientations viewed as part of $c_1^t $ or $c_2^t
$). The simplices of $\partial \zeta^t = $ again decompose into those
contained in $\Delta \setmin D $ for all $t $ and those running into
$D $. The first type will be denoted by $\gamma^t $, the second type
by $\nu^t $. Applying the usual Stokes formula to the integrals over
$c_i^t $ gives
$$
\Int_{c_1^t} < \nabla \eta_i, \ve > - \Int_{\partial c_i^t-D} <
\eta_1, \ve > = \Int_{\nu_i^t} < \eta_i, \ve > + (-1)^i
\Int_{\zeta^t} < \eta_i, \ve > \ ,
$$
and $\alpha = \lim_{t \to 0} ( \int_{\nu_1^t} \!\!< \!\! \eta_1, \ve
\!\!>\!\! - \int_{\nu_2^t} \!\!<\!\! \eta_2, \ve \!\!>\!\! +
\int_{\zeta^t} \!\!<\!\! \eta_1 + \eta_2, \ve \!\!>\!\! +
\int_{\gamma^t} \!\!<\!\! \omega, \ve \!\!>\!\! ) $.
Now, $\eta_1 + \eta_2
= \nabla \omega $ and the Stokes formula gives
\begin{equation} \label{eq:alphastokes} \alpha = \lim_{t \to 0} (
  \Int_{\nu_1^t} < \eta_1, \ve > - \Int_{\nu_2^t} < \eta_2, \ve >
  + \Int_{\nu^t} <\omega, \ve >) \ .
\end{equation}
With this presentation of $\alpha $ independence of the choices made
is easily shown. If e.g. $\nabla \omega= \eta_1 + \eta_2 = \eta'_1 +
\eta'_2 $ are two different choices, it follows that $\eta_1 - \eta'_1
= \eta'_2 - \eta_2 \in L_1^p \cap L_2^p = \Emero^p $, and hence the
rapid decay of $\ve $ implies that $\int_{\nu_i^t} < \eta_i -
\eta'_i, \ve > = 0 $.

If $\omega \equiv 0 \, {\rm mod} \, L_1+L_2 $, say $\omega =
\omega_1 + \omega_2 \in L_1^{p-1} + L_2^{p-1} $, one can take
$\eta_i := \nabla \omega_i $ and obtains that $\alpha=0 $, since
$\partial \nu_i^t $ consists of $\nu^t $ and some simplices running
into $D \setmin 0 $, let us denote them by $\mu^t $. Then
$$
\lim_{t \to 0} \Int_{\nu_1^t} < \nabla \omega_1, \ve > = \lim_{t \to
  0} (\Int_{\nu^t} < \omega_1, \ve > + \Int_{\mu^t} < \omega_1, \ve
>) = \lim_{t \to 0} \Int_{\nu^t} < \omega_1, \ve > \ ,
$$
the last equality following from the rapid decay of $\ve $ and the at
worst meromorphic behavior of $\omega_1 $ in the limit process
involved. Taking care of the orientations (that are all induced by the
one of $c $) gives that $\alpha=0 $ in this case. The cases $\omega =
\nabla \Omega $ or the corresponding cases to show independence of the
choices of the rapid decay chains are similar and omitted here.
\\ \qed 

\subsubsection{Compatibility with the long exact sequences}

We now want to show that the pairings defined in the last section
commute with the mappings of the long exact sequences \eqref{eq:les1}
and \eqref{eq:les2}. The first step toward this has already been done
in Lemma \ref{lemma:comp1}.

First, we take a more precise look at the rd-homology sequence
\eqref{eq:les2}. In it we find the following mapping, the first row of
the diagram:
\begin{equation} \label{eq:v1v2delta}
  \begin{array}{ccccc}
    \Hrd_{p+1}(V_1, \partial V_1; D) & \!\! \oplus \!\! & \Hrd_{p+1}(V_2,
    \partial V_2; D) & \lto & \Hrd_{p+1}(\Delta, \partial \Delta; D) \\
    \uparrow \cong & & \uparrow \cong & & \uparrow \cong\\
    \Hrd_{p+1}(\Delta, \Delta \setmin D_1; D_1) & \!\! \oplus \!\! & \Hrd_{p+1}(\Delta,
    \Delta\setmin D_2; D_2) & \lto & \Hrd_{p+1}(\Delta, \Delta \setmin 0; D)\\[0.1cm]
    [c_1 \otimes \ve_1] & \!\! + \!\! & [c_2 \otimes \ve_2] & \mapsto &
    [c_1 \otimes \ve_1] - [c_2 \otimes \ve_2]
  \end{array} \ .
\end{equation}
The vertical mapping at the right hand side is defined by 'capping off
at $\partial \Delta $'. More precisely, let $c \otimes \ve \in
\Crd_{p+1}(\Delta) $. By subdivision, we can decompose the topological
chain $c $ into a sum $c = \widetilde{c} + \gamma $ such that
$\widetilde{c} = c \cap \Delta $ and we have $\Hrd_{p+1}(\Delta,
\Delta \setmin 0; D) \stackrel{\cong}{\to} \Hrd_{p+1}(\Delta, \partial
\Delta; D) $, $[c \otimes \ve] \mapsto [\widetilde{c} \otimes \ve] $.
Now suppose we have a $\omega \in \Eess^{p-1} $ with $\nabla \omega =
\eta_1 + \eta_2 \in L_1^p+L_2^p $ (and therefore also $\nabla \eta_1
= - \nabla \eta_2 \in \Emero^p $) and $c_\nu \otimes \ve_\nu \in
\Crd_{p+1}(V_\nu) $ with $\partial(c_\nu \otimes \ve_\nu) \in
\Crd_p(\partial V_\nu) + \Crd_p(D) $.

We decompose $c_\nu = \widetilde{c}_\nu + \gamma_\nu $ as before. Then
$\gamma_\nu \cap D_{\neq \nu} = \emptyset $ and the 'limit Stokes
formula'
$$
\Int_{\gamma_1} < \ve_1, \nabla \eta_1 > = - \Int_{\gamma_1} <
\ve_1, \nabla \eta_2 > = - \Int_{\partial \gamma_1- D_1} < \ve_1,
\eta_2 >
$$
tells us that
\begin{equation} \label{eq:comp2} \Int_{c_\nu} < \ve_\nu, \nabla
  \eta_\nu > - \Int_{\partial c_\nu- D_\nu} < \ve_\nu, \eta_\nu > =
  \Int_{\widetilde{c}_\nu \phantom{tld}} < \ve_\nu, \nabla \eta_\nu
  > - \Int_{\partial \widetilde{c}_\nu- D_\nu} < \ve_\nu, \eta_\nu
  > \ .
\end{equation}

\begin{lemma} \label{lemma:comp2} The maps in the diagram
  $$
\begin{array}{ccc}
  \Hdr^p(\mod{L_1}{\Emero}) \oplus \Hdr^p(\mod{L_2}{\Emero}) & \lot &
  \Hdr^{p-1}(\mod{\Eess}{L_1+L_2}) \\[0.5 em] 
  \Hrd_{p+1}(V_1, \partial V_1; D_1) \oplus \Hrd_{p+1}(V_2, \partial V_2; D_2)
  & \lto & \Hrd_{p+1}(\Delta, \partial \Delta; D)
\end{array}
$$
are compatible with the given pairings.
\end{lemma}
\pf With the notation introduced right above the lemma, the diagram
maps the given elements as $[\eta_1] + [\eta_2] \ot [\omega] $ and
$[c_1 \otimes \ve_1] + [c_2 \otimes \ve_2] \mapsto [\widetilde{c}_1
\otimes \ve_1] - [\widetilde{c}_2 \otimes \ve_2] $. The desired
equation $< [\omega], [\widetilde{c}_1 \otimes \ve_1] > - <
[\omega], [\widetilde{c}_2 \otimes \ve_2] > = < [\eta_1], [c_1
\otimes \ve_1] > + < [\eta_2], [c_2 \otimes \ve_2] > $ follows
directly from the
definitions using \eqref{eq:comp2}.
\\ \qed 

\noindent
It remains to prove
\begin{lemma}
  \label{lemma:comp3} The maps in the diagram
  $$
\begin{array}{ccc}
  \Hdr^{p-1}(\mod{\Eess}{L_1+L_2}) & \lot &
  \Hdr^{p-1}(\mod{\Eess}{\Emero}) \\[0.5em]
  \Hrd_{p+1}(\Delta, \partial \Delta; D) & \lto & \Hrd_p(\Delta, \Delta
  \setmin D; D)
\end{array}
$$
are compatible with the given pairings.
\end{lemma}
\pf We compute the connecting morphism $\psi $ in the long exact
sequence induced by
$$
0 \to \Crd_\ast(\Delta, \Delta \setmin D; D) \to \bigoplus_{i=1,2}
\Crd_\ast(\Delta, \Delta \setmin D_i; D) 
\to \Crd_\ast(\Delta, \Delta \setmin 0; D) \to 0 \ .
$$
If we write $(\cc_\nu)_p := \Crd_p(\Delta \setmin D_\nu) + \Crd_p(D) $
and $\cc_p := \Crd_p(\Delta) $, we have $\Crd_p(\Delta \setmin 0) +
\Crd_p(D) = (\cc_1)_p + (\cc_2)_p $ and $\Crd_p(\Delta \setmin D) +
\Crd_p(D) = (\cc_1)_p \cap (\cc_2)_p $. Therefore the homological yoga
to determine $\psi([c \otimes \ve]) $ for a given chain $c \otimes \ve
\in \cc_p $ with $\partial (c \otimes \ve) \in (\cc_1 + \cc_2)_{p-1} $
reads:
$$
\begin{array}{rccccccl}
  &&& (\mod{\cc}{\cc_1})_p \oplus (\mod{\cc}{\cc_2})_p & \to &
  (\mod{\cc}{\cc_1 + \cc_2})_p & \to 0\\ [0.1cm]
  &&& \downarrow && \downarrow\\[0.1cm]
  0 \to & (\mod{\cc}{\cc_1 \cap \cc_2})_{p-1} & \to & (\mod{\cc}{\cc_1})_{p-1}
  \oplus (\mod{\cc}{\cc_2})_{p-1} & \to & (\mod{\cc}{\cc_1 +
    \cc_2})_{p-1} & \to 0
\end{array}
$$
where $c \otimes \ve $ makes its way through this diagram as follows
(to be read from the upper right to the lower left corner):
\begin{equation} \label{eq:diagchase}
\begin{array}{rccccccl}
  &&& 0 + [-c \otimes \ve] & \mapsto &
  [c \otimes \ve] \\ [0.1cm]
  &&& \downarrow && \downarrow\\[0.1cm]
  0 \mapsto & [- \partial (c_1 \otimes \ve)] & \mapsto & [-\partial (c_1
  \otimes \ve)] + [-\partial (c_1 \otimes \ve)] & \mapsto & 0 \ .
\end{array}
\end{equation}
Here we decomposed the topological chain $c \in C_p(\Delta) $ into the
sum $c=c_1 + c_2 $ with $\partial c_\nu \subset \Delta \setmin D_\nu
$. It follows, that $0 + [-\partial (c \otimes \ve)] = [-\partial (c_1 \otimes
\ve)] + [-\partial (c_1 \otimes \ve)] \in (\mod{\cc}{\cc_1})_{p-1}
\oplus (\mod{\cc}{\cc_2})_{p-1} $.

Now, let $\omega \in \Eess^{p-1} $ be given such that $\nabla \omega
\in \Emero^p $. To keep notation, we let $\eta_1 := \nabla \omega \in
\Emero^p \subset L_1^p $. Suppose we started in \eqref{eq:diagchase}
above with $c \otimes \ve \in \Crd_p(\Delta) $ such that $\partial c
\subset \Crd_{p-1}(\partial \Delta) + \Crd_{p-1}(D) $.  Then we have
to consider the diagram
$$
\begin{array}{ccccccc}
  \Hdr^{p-1}(\mod{\Eess}{L_1+L_2}) & \ot & \Hdr^{p-1}(\mod{\Eess}{\Emero}) &
  \phantom{shift} & [\omega] & \ot & [\omega] \\
  \times && \times \\
  \Hrd_{p+1}(\Delta, \partial \Delta; D) & \to & \Hrd_p(\Delta, \Delta
  \setmin D; D) & &
  [c \otimes \ve] & \mapsto & [-\partial (c_1 \otimes \ve)] \ .
\end{array}
$$
Starting with the left hand side we obtain 
$$
< [\omega], [c \otimes
\ve] > \stackrel{\mbox{\tiny def}}{=} 
\alpha(\eta_1, 0, c_1 \otimes \ve, c_2 \otimes \ve) = - \!\!\!\!\!\!
\Int_{\partial c_1 \cap \partial \Delta - D} \!\!\!\!\!\!< \ve,
\eta_1 > + \!\!\!\!\!\! \Int_{\partial(\partial c_1 \cap \partial
  c_2) - D} \!\!\!\!\!\!\!\!\!\!< \ve, \omega > = < [\omega], [-
\partial (c_1 \otimes \ve)] > \ ,
$$
which is the right hand side.
\\ \qed 

\subsection{The Localization Lemma}

In the section above, we have constructed a pairing between the long
exact sequences \eqref{eq:les1} and \eqref{eq:les2}. It follows, that
in order to prove perfectness of the irregularity pairing, it suffices
to do so for $\Hdr^p(\mod{\Eess}{L_1 + L_2}) \times
\Hrd_{p+1}(\Delta, \partial \Delta; D) \to \C $ and
$\Hdr^{p-1}(\mod{\Enu}{\Emero}) \times \Hrd_{p+1}(V_\nu, \partial
V_\nu; D_\nu) \to \C $.

\noindent
Now, consider the short exact sequence of de Rham complexes
\begin{equation} \label{eq:ses1} 0 \ot \deR{\mod{\Eess}{L_1 +
      L_2}} \ot \deR{\mod{\Eess}{L_2}} \ot
  \deR{\mod{L_1}{\Emero}} \ot 0
\end{equation}
where the mapping is defined via the isomorphism
$\mod{L_1^p}{\Emero^p} \stackrel{\cong}{\to} \mod{L_1^p +
  L_2^p}{L_2^p} $, which follows from $L_1^p \cap L_2^p= \Emero^p
$. In the same manner we have a short exact sequence of complexes of
rd-chains:
\begin{equation} \label{eq:ses2} 0 \to \Crd_\ast(\Delta \setmin D_2,
  \Delta \setmin D; D) \to \Crd_\ast(\Delta, \Delta \setmin D_1; D)
  \to \Crd_\ast(\Delta, \Delta \setmin 0; D) \to 0 \ .
\end{equation}
Excision and retraction to the boundary again induces a natural isomorphism
$\Hrd_p(\Delta \setmin D_2, \Delta \setmin D; D)
\cong \Hrd_p(V_1 \setmin D_2,
\partial V_1; D) $. Again, for given
\begin{enumerate}
\item $\omega \in \Eess^p $ with $\nabla \omega \in L_2^p $ and
\item $c \otimes \ve \in \Crd_p(V_1 \setmin D_2) $ with $\partial c
  \in C_{p-1}(\partial V_1 \setmin D_2) + C_{p-1}(D) $.
\end{enumerate}
we define
$$
< [\omega], [c \otimes \ve] > := \Int_c < \ve, \nabla \omega > -
\Int_{\partial c_1- D_1} < \ve, \omega > \ ,
$$
where the first integral exists, because $c $ does not intersect $D_2
$ and $\nabla \omega $ is meromorphic along $D_1 $, and the second one
does, as $\partial c - D_1 $ lies in $X\setmin D $. With the same
arguments already used several times above one shows that this gives
indeed a well-defined pairing of the long exact sequences
corresponding to \eqref{eq:ses1} and \eqref{eq:ses2}.

\bigskip
We have thus reduced the problem of perfectness of the period pairing
at a crossing-point to the perfectness of a pairing involving
(co-)homology groups supported at the crossing-point $0 $ and a
pairing involving groups supported on the complement $D \setmin \{ 0
\} $. Together with the considerations at the beginning of this
section, we have now achieved the Localization Lemma, which we will
state using the following notion:
\begin{definition} \label{def:locperf} We say that the rank $1 $
  connection $(L, \nabla) $ on $X $, singular along $D $, satisfies
  {\bf local perfectness}, if there is an open covering ${\mathscr U}
  $ of the tubular neighborhood $V $ of $D $ consisting of small
  enough bi-discs $\Delta $, such that the following holds:
  \begin{enumerate}
  \item for each local situation $D= D_1= \{ x_1=0 \} \subset \Delta
    $, the pairing
$$
\Hdr^{p-1}(\Delta; \mod{\Eess}{\Emero}) \times \Hrd_p(\Delta, \Delta
\setmin D_1; D_1) \lto \C \mbox{\quad is perfect, and}
$$
\item for each local situation $D=\{ x_1x_2=0 \} \subset \Delta $, the
  pairings
    $$
    \Hdr^{p-1}(\Delta; \mod{\Eess}{\Enu}) \times \Hrd_p(\Delta \setmin
    D_{\neq \nu}, \Delta \setmin D; D) \to \C \ ,
$$
$(\nu=1,2) $, and
$$ \Hdr^{p-1}(\Delta; \mod{\Eess}{L_1 + L_2}) \times
\Hrd_{p+1}(\Delta, \Delta \setmin D; D) \to \C \ ,
$$
are perfect.
\end{enumerate}
\end{definition}
\begin{lemma} \label{thm:local} Let $X $ be of dimension
  $\dim(X)=2 $ and $(L, \nabla) $ be an integrable meromorphic rank $1
  $ connection on $X $ with singularities along the normal crossing
  divisor $D $. Then the period pairing $\Hdr^p(U; L(\sD),\nabla) \times
  \Hrd_p(X; L^\vee) \to \C $ is perfect if and only if $(L, \nabla) $
  satisfies local perfectness.
\end{lemma}

\section{Local perfectness} \label{sec:localperf}

In any of the
local situations above, consider the
formal connection $\widehat{L} := L \otimes \O_{\widehat{X|Y}}(\sD) $,
where $\O_{\widehat{X|Y}} $ denotes the formal completion of $\O_X $
with respect to the stratum $Y $ of $D $ considered, namely $Y $ being
a smooth component of $D $ or a crossing-point. By assumption
(Definition \ref{def:good}), $\widehat{L}|_\Delta $ is isomorphic to
the formal completion of a connection of the form $e^\alpha \otimes R
$, with $\alpha(x):= x_1^{-m_1} x_2^{-m_2} u(x) $ such that $u(0) \neq
0 $ and $R $ a regular singular connection. It follows that $(L,
\nabla)|_\Delta $ itself is isomorphic to $e^\alpha \otimes R $ (since
$L \otimes e^{-\alpha} $ is necessarily regular singular). 

\subsection{Local pairing supported on a crossing-point}
\label{sec:locperfcross}

We consider the local situation ii) of Definition \ref{def:locperf},
i.e. $\Delta $ is a small bi-disc around the crossing-point $0 \in D
$.

\subsubsection{Local rd-homology}

We introduce the following notions (after chosing fixed local
coordinates) in the situation of a connection of the form $e^\alpha
\otimes R $ with $\alpha(x):= x_1^{-m_1} x_2^{-m_2} u(x) $ with $u(0)
\neq 0 $ and $R $ a regular singular connection.
\begin{definition} \label{def:stokes} The {\bf Stokes directions} of
  $e^\alpha \otimes R $ at $0 $ are the elements of
  $$
  \Sigma_0 := \{ (\vt_1, \vt_2) \mid -m_1 \vt_1 -m_2 \vt_2 + {\rm
    arg}(u(0)) \in (\frac{\pi}{2}, \frac{3 \pi}{2}) \} \subset
  \pi^{-1}(0) \cong S^1 \times S^1 \ .
  $$
  The {\bf Stokes bisectors} are the bisectors $\S_0 :=
  \bigcup_{(\vt_1, \vt_2) \in \Sigma_0} \R^+ e^{i\vt_1} \times \R^+
  e^{i\vt_2} \subset \C^2 $.
\end{definition}

\noindent
We now calculate the local rapid decay homology groups at a crossing
point. In the following, we will abbreviate $\Hrd_\ast(Y, \partial) :=
\Hrd_\ast(Y, \partial Y; L) $ and in the same way for the rapid decay
chains $\Crd_\ast(Y, \partial) $ for any $Y \subset \Delta $.
\begin{proposition} \label{prop:rdcrossiso} In the situation above,
  we have
  $$
  \Hrd_\ast(\Delta, \partial) \cong H_\ast(\Delta \cap \S_0,
  \partial \Delta \cup D) \ ,
$$
where the right hand side denotes the usual singular homology group.
\end{proposition}
\pf We first observe, that we have a natural homomorphism
\begin{equation} \label{eq:morph}
H_\ast(Y \cap \S_0, \partial Y \cup D) \lto
\Hrd_\ast(Y,\partial)
\end{equation}
for any $Y \subset \Delta $, since $e^\alpha $ is rapidly decaying
along any simplex contained in the Stokes region. We claim, that this
is an isomorphism for $Y=\Delta $.

Now, both sides of the morphism \eqref{eq:morph} can be embedded in
Mayer-Vietoris sequences by decomposing $\Delta $ into bisectors in
as follows. Let $\Delta_1 = \bigcup_{i=1}^n S(\nu_i) $ and
$\Delta_2 = \bigcup_{j=1}^m S(\mu_j) $ be a decomposition of the disc
into closed sectors $S(\nu_i) $ and $S(\mu_j) $, where the $\nu_i $
and $\mu_j $ are intervals in $S^1 $, cyclically ordered, such that
$\bigcup_{i=1}^n \nu_i = S^1 = \bigcup_{j=1}^m \mu_j $. We denote by
$$
[0, p_{i+1}] := \l_{i+1} := S(\nu_i) \cap S(\nu_{i+1}) \mbox{\quad and
  \quad} [0, q_{j+1}] := \k_{j+1} := S(\mu_j) \cap S(\mu_{j+1})
$$
the common lines of two consecutive sectors. We assume that the
intersection $(\nu_i \times \mu_j)
\cap \S_0 $ is either empty or connected. The decomposition $
\Delta = \bigcup_{j=1}^m \Delta_1 \times S(\mu_j) =: \bigcup_{j=1}^m
Z_j $ gives rise to the following Mayer-Vietoris sequence with $Z_{ij}
:= Z_i \cap Z_j $:
$$
0 \ot \Crd_\ast(\Delta, \partial) \ot \bigoplus_{j=1}^m \Crd_\ast(Z_i,
\partial) \ot \bigoplus_{1 \le i< j \le m} \Crd_\ast(Z_{ij}, \partial)
\ot \ldots \ .
$$
Because of $Z_{ijk} = \Delta_1 \times (S(\mu_i) \cap S(\mu_j) \cap
S(\mu_k)) = \Delta_1 \times 0 \subset D $, this sequences reduces to a
short exact sequence (recall that in the definition of the rapid
decaying chains, we have moduled out all chains contained in $D $) and
we obtain the following long exact sequence of rapid decay homology
groups:
\begin{equation} \label{equ:MV1} \ldots \to \Hrd_{n+1}(\Delta,
  \partial) \to \bigoplus_{1 \le i < j \le m} \Hrd_n(Z_{ij}, \partial)
  \to \bigoplus_{i=1}^m \Hrd_n(Z_i,
  \partial)   \to   \Hrd_n(\Delta, \partial)   \to   \ldots 
\end{equation}
The same holds for the left hand side of \eqref{eq:morph} and thus
\eqref{eq:morph} induces a natural morphism between these long exact
sequences.

In a similar way, we fix $j $ and consider the decomposition of
$\Delta_1 \times S(\mu) := Z := Z_j $ according to the decomposition
of $\Delta_1 $ from above:
$$
Z = \Delta_1 \times S(\mu) = \bigcup_{i=1}^n S(\nu_i) \times S(\mu) \
.
$$
Let $B_i := S(\nu_i) \times S(\mu) $. The analogous Mayer-Vietoris
sequence of the $\Crd $-groups is again short exact (as $B_{ijk} :=
B_i \cap B_j \cap B_k \subset D $) and therefore gives rise to the
following exact Mayer-Vietoris sequence:
\begin{equation} \label{equ:MV2} \ldots \to \Hrd_{n+1}(Z, \partial)
  \to \bigoplus_{1 \le i < k \le n} \Hrd_n(B_{ik}, \partial) \to
  \bigoplus_{i=1}^m \Hrd_n(B_i,
  \partial) \to \Hrd_n(Z, \partial) \to \ldots \ ,
\end{equation}
together with a map form the corresponding MV-sequence for $H_\ast(Z
\cap \S_0, \partial Z \cup D) $.  Observe that
$$
B_{ik} = \left\{
  \begin{array}{ll}
    0 \times S(\mu) \subset D & \mbox{for } |i-k| \ge 2 \\
    \l_{i+1} \times S(\mu) & \mbox{for } k=i+1
  \end{array} \right. \ \mbox{ and } 
Z_{ij} = \left\{
  \begin{array}{ll}
    \Delta_1 \times 0 \subset D & \mbox{for } |i-j| \ge 2 \\
    \Delta_1 \times \k_i  & \mbox{for } j=i+1
  \end{array} \right. \ .
$$
We keep $j $ fixed and decompose $\Delta_1 \times \k = \bigcup_{i=1}^n
S(\nu_i) \times \k $ for $\k = \k_j $. This leads to the
Mayer-Vietoris sequence
\begin{equation} \label{equ:MV3} \ldots \to \Hrd_{n+1}(Z_{i,i+1},
  \partial) \to \bigoplus_{1 \le k < l \le n} \Hrd_n(C_{kl}, \partial,
  D) \to \bigoplus_{k=1}^m \Hrd_n(C_k,
  \partial) \to \ldots
\end{equation}
where $C_k := S(\nu_k) \times \k $ and thus $ C_{kl} := C_k \cap C_l =
\left\{
  \begin{array}{ll}
    0 \times \k \subset D & \mbox{for } |k-l| \ge 2 \\
    \l_k \times \k  & \mbox{for } l=k+1
  \end{array} \right. 
$. The similar assertion holds for the left hand side of
\eqref{eq:morph}. The proposition now follows from the following
lemma:  
\begin{lemma} \label{lemma:lk} There are natural isomorphisms
  $$
  \begin{array}{llll} {\rm i)} & \Hrd_\ast(\l \times \k, \partial) &
    \to & H_\ast((\l \times \k) \cap
    \S_0, \partial (\l \times \k) \cup D) \mbox{\quad and}\\[0.1cm]
    {\rm ii)} & \Hrd_\ast(S(\nu) \times \k, \partial) & \to &
    H_\ast((S(\nu) \times \k) \cap \S_0, \partial (S(\nu)
    \times \k) \cup D) \ .
  \end{array}
$$
\end{lemma}
\pf Recall that $(L, \nabla) \cong
  e^\alpha \otimes R $ for some $\alpha=x_1^{-m_1} x_2^{-m_2} u(x) $ with
$u(0,0) \neq 0 $ . We start with the proof for i). We
claim, that if $\l \times \k \not\subset \S_0 $, then there is
no rapidly decaying chain $c \otimes \ve $ approaching $0 $ inside $\l
\times \k $, which is not entirely contained in $D $. For, let
$l=[0,e^{i\vt_1}] $ and $\k=[0,e^{i\vt_2}] $ and suppose that $\ve =
e^{\alpha(x)} $ is rapidly decaying as $x $ varies in $c $, it follows
that
$$
\exp(|x_1|^{-m_1} |x_2|^{-m_2} \cdot |u(x)| \cdot \cos(-m_1 \vt_1 -m_2
\vt_2 + {\rm arg}(u(x)))) \le C_{NM} |x_1|^N |x_2|^M
$$
for $x \in c $. But then $-m_1 \vt_1 -m_2 \vt_2 + {\rm arg}(u(x)) \in
(\pi/2, 3 \pi/2) $ for $x $ small enough and therefore $(\vt_1, \vt_2)
\in \overline{\Sigma}_0 $ and, because of the properties of the chosen
decomposition with respect to the Stokes bisectors, finally $(\vt_1,
\vt_2) \in \Sigma_0 $, thus $\l \times \k \subset \S_0 $. For $\l
\times \k \subset \S_0 $, however, the assertion is clear.

In order to prove part ii), consider $(S(\nu) \times \k) \cap
\S_0 $, which is either empty or connected by assumption. Let
$\nu =: [\xi, \xi'] \subset S^1 $ and let $\zeta $ be the direction of
the line $\k $, i.e. $\k=[0, e^{i \zeta}] \subset \C $. Then $(S(\nu)
\times \k) \cap \S_0 $ is the union of the radii with
directions contained in $(\nu \times \{ \zeta \}) \cap \Sigma_0
\subset S^1 \times S^1 $, with $\Sigma_0 $ being the Stokes
directions. Let $\rho $ be the interval $\rho := \nu \cap
\Sigma_0 $. If $\rho $ is empty or $\rho=\nu $ the assertion is
clear. Otherwise, assume that $\rho=[\xi, \eta] \subsetneq [\xi, \xi']
= \nu $. Let $\l $ be the radius with direction $\xi $, so that we now
have $\l \times \k \subset \S_0 $. If $h: [0,1] \times (\nu
\times \{ \zeta \}) \to (\nu \times \{ \zeta \}) $ denotes the linear
retraction of $\nu $ to $\{ \zeta \} $, i.e.
$h(t,(x,\zeta))=((1-t)x+t\xi, \zeta) $, then
$$
H: [0,1] \times (S(\nu) \times \k) \to (S(\nu) \times \k) \ , \ (t,
(r_1 e^{ix}, r_2 e^{i\zeta})) \mapsto (r_1 e^{ih(t,x,\zeta)}, r_2
e^{i\zeta})
$$
retracts $S(\nu) \times \k $ to $\l \times \k \subset \S_0 $.
We claim, that $H $ preserves the rapid decay condition. To prove
this, consider a curve $\gamma: [0,1] \to S(\nu) \times \k $ with
$\gamma^{-1}(D)=0 $, such that $\ve=e^{\alpha(x)} $ is rapidly
decaying along $\gamma $. We have to show that $\ve $ then is rapidly
decaying along all the curves $H(t, .) \circ \gamma $, $t \in [0,1] $.
Let $\gamma(s)=(r_1(s) e^{i \theta(s)}, r_2(s) e^{i \zeta}) $. From
the rapid decay of $\ve $ along $\gamma $, it is clear that the
direction of $\gamma $ at the point $\gamma(0) \in D $ is in the
closure of the Stokes-directions, i.e. $(\theta(0), \zeta) \in
\overline{\Sigma_0} $. But then, the direction of the transposed
curve $H(t,.) \circ \gamma $ at the point $H(t, \gamma(0)) $ by
construction lies in $\Sigma_0 $, whence $\ve $ is rapidly
decaying along this curve. Thus, $H $ induces isomorphisms
$$
\Hrd_\ast(S(\nu) \times \k, \partial) \cong \Hrd_\ast(\l \times \k,
\partial) \cong H_\ast((\l \times \k) \cap \S_0, \partial)
\cong H_\ast((S(\nu) \times \k) \cap \S_0,
\partial) \ .
$$
\\ \qed 

The proposition now follows immediately using the 5-lemma applied to
the morphism \eqref{eq:morph} induces on the various Mayer-Vietoris
sequences above (\eqref{equ:MV1} -- \eqref{equ:MV3}).
\\ \qed

\medskip\noindent Computing the singular homology groups appearing in
the proposition, we obtain the following
\begin{theorem} \label{thm:rddim} For a rank $1 $ connection $L $
  with formal model $e^\alpha \otimes R $, we have
  $$
  \dim \Hrd_\ast(\Delta, \partial) = \left\{
    \begin{array}{lll}
      0 & for & \ast=1 \mbox{ or } \ast \ge 4 \\[0.3cm]
      (m_1, m_2) & 
      \mbox{for} & \ast=2 \mbox{ or } 3, 
    \end{array} \right.
$$
where $\alpha=x_1^{-m_1} x_2^{-m_2} \cdot u_\alpha(x) $ with
$u_\alpha(0) \neq 0 $ and $(m,n) $ denotes the greatest common divisor
of two non-negative integers $m,n $. Note, that the homology groups in
degree $0 $ will play no role here.
\end{theorem}
\pf The subspace of the Stokes directions $\Sigma_0 \subset S^1 \times
S^1 $ of $\alpha $ are homotopy equivalent to a torus knot of type
$(m_1, m_2) $, which we denote by ${\cal K} $. Let ${\cal R} \subset
\Delta $ be the union of the radial sheets with directions in ${\cal
  K} $, a radial sheet being the product of the two radii in each
direction, then the homology group to be computed is isomorphic to
$H_\ast({\cal R}, \partial \Delta \cup D) $. Now, consider a
decomposition $\Delta= \bigcup S(\nu_i) \times S(\mu_j) $ in bisectors
as above, where we assume that the intersection of $\S_0 $ with any
bisector has at most one connected component. Now observe that for two
radii $\l $ and $\k $ and a sector $S(\nu) $ one has homotopy
equivalences
\begin{equation} \label{equ:S2} \mod{(\l \times \k)}{\partial \Delta
    \cup D} \simeq S^2 \mbox{\quad and \quad} \mod{(l \times
    S(\nu))}{\partial \Delta \cup D} \simeq S^2 \ .
\end{equation}
Starting with this observation, we can make our way through the
various Mayer-Vietoris sequences induced from the decomposition of
$\Delta $ in direct analogy to the sequences \eqref{equ:MV1} --
\eqref{equ:MV3}. One easily deduces that the homology groups to be
computed vanish in degree $1 $ and greater than or equal to $4 $. As
for the remaining degrees $2,3 $, one sees that the one-dimensional
contributions coming from $H_2(S^2)=\C $ via \eqref{equ:S2}
distinguish each other as long as they come from points on the torus
knot ${\cal K} $ (i.e. the endpoint of $\l \times \k $ or any
direction of some $S(\nu) \times \k $ respectively) that belong to the
same connected component of the torus knot. A careful book-keeping
thus gives the desired result for $H_2 $ and $H_3 $, the integers
$(m_1, m_2) $ being the number of connected components of the torus
knot. (Another way to look at it, is to cut the torus along the knot
${\cal K} $ and to use a Mayer-Vietoris argument to see that the
dimension of the homology groups in question are given as the
dimension of the analogous homology groups in the special case
$m_1=1=m_2 $, which is easily seen to be $1 $, times the number of
connected components of the torus knot
${\cal K} $, namely $(m_1, m_2) $).
\\ \qed 

\subsubsection{Local de Rham cohomology}

We are now going to compute the local de Rham cohomology
$H^\ast_{dR}(\Eess/L_1+L_2) $ at the crossing point, where the rank
one connection $L $ has the formal model $e^\alpha \otimes R $
with $\alpha = x_1^{-m_1} x_2^{-m_2} \cdot u_\alpha(x_1,
x_2) $ with $u_\alpha(0,0) \neq 0 $.

Since the sheaf $\Eess/L_1+L_2 $ has support at the origin, it is
sufficient to consider the stalk at $0 $. We will keep the same
notation $\Eess, L_1 $ and $L_2 $ but in the following think of
these as the stalks. The main result of this section is the following
\begin{theorem} \label{thm:dRdim} In the situation above, we have
  $$
  \dim \Hdr^\ast(\mod{\Eess}{L_1+L_2}) = \left\{
    \begin{array}{lll}
      0 & \mbox{for} & \ast \neq 0,1 \\[0.3cm]
      (m_1, m_2) & \mbox{for} & \ast = 0,1
    \end{array} \right.
$$
\end{theorem}
\pf The regular singular part does not give any contribution to
$\mod{\Eess}{L_1+L_2} $, so that we can omit it. Since $u_\alpha(0) \neq 0
$, we can locally transform the connection by multiplication with
$u_\alpha^{-1} $ and obtain a new connection $L' $ whose induced complex
$\mod{\Eess'}{L_1'+L_2'} $ is quasi-isomorphic to the analogous
complex for $e^\alpha $. Therefore we can also omit the factor $u_\alpha(x) $
and take $\alpha= x_1^{-m_1} x_2^{-m_2} $. We then have to consider
the following complex
$$
\mod{\Eess}{L_1+L_2} \stackrel{\vi}{\lto} \mod{\Eess}{L_1+L_2}
\oplus \mod{\Eess}{L_1+L_2} \stackrel{\psi}{\lto}
\mod{\Eess}{L_1+L_2}
$$
with
\begin{align*}
  & \vi(u) := (\frac{\partial u}{\partial x_1} - m_1 x_1^{-m_1-1}
  x_2^{-m_2} u \, , \, \frac{\partial u}{\partial x_2} - m_2
  x_1^{-m_1} x_2^{-m_2-1} u) \mbox{\quad and} \\
  & \psi(\omega_1, \omega_2 ) := \frac{\partial \omega_1}{\partial
    x_1} - \frac{\partial \omega_2}{\partial x_2} + m_2 x_1^{-m_1}
  x_2^{-m_2-1} \omega_1 - m_1 x_1^{-m_1-1} x_2^{-m_2} \omega_2 \ .
\end{align*}

\medskip\noindent{\em 1. Step:} If we let
$$
K_{MN} := \{ \sum_{k,l \in \Z} f_{kl} x_1^k x_2^l \in L_1+L_2 \mid
f_{kl}=0 \mbox{ for } k<M \mbox{ and } l<N \} \subset L_1+L_2 \ ,
$$
it follows that $K_{1,1} \stackrel{\vi}{\lto} K_{-m_1, -m_2+1} \oplus
K_{-m_1+1, -m_2} \stackrel{\psi}{\lto} K_{-2m_1, -2m_2} $.  Consider
the following diagram with exact rows:
\begin{equation} \label{eq:Kdiag}
  \begin{array}{rcccccccl}
    && 0 && 0 && 0 \\
    && \downarrow && \downarrow && \downarrow\\
    0 & \to & \mod{L_1+L_2}{K_{1,1}} & \to & \mod{\Eess}{K_{1,1}} & \to &
    \mod{\Eess}{L_1+L_2} & \to & 0 \\
    && \downarrow && \downarrow && \downarrow \\
    0 & \to & 
    \footnotesize \begin{gathered} \mod{L_1+L_2}{K_{-m_1, -m_2+1}} \qquad
      \\[-.1cm] \oplus \\[-.1cm] \qquad \mod{L_1+L_2}{K_{-m_1+1, -m_2}}
    \end{gathered} \normalsize & \to & \footnotesize \begin{gathered}
      \mod{\Eess}{K_{-m_1, -m_2+1}} \qquad
      \\[-.1cm] \oplus\\[-.1cm]
      \qquad \mod{\Eess}{K_{-m_1+1, -m_2}} \end{gathered} \normalsize
    & \to & \footnotesize \begin{gathered} \mod{\Eess}{L_1+L_2}
      \qquad
      \\[-.1cm] \oplus\\[-.1cm]
      \qquad \mod{\Eess}{L_1+L_2} \end{gathered} \normalsize
    & \to & 0 \\
    && \downarrow && \downarrow && \downarrow \\
    0 & \to & \mod{L_1+L_2}{K_{-2m_1, -2m_2}} & \to &
    \mod{\Eess}{K_{-2m_1,-2m_2}} & \to &
    \mod{\Eess}{L_1+L_2} & \to & 0 \\
    && \downarrow && \downarrow && \downarrow\\
    && 0 && 0 && 0
  \end{array}
\end{equation}
where the vertical arrows are induced by the connection $\nabla $. We
denote the columns of this diagram as ${\cal K}' $, ${\cal K} $ and
$\overline{{\cal K}} $, so that the diagram reads as the short exact
sequence $0 \to {\cal K}' \to {\cal K} \to \overline{{\cal K}} \to 0 $.

\begin{lemma} \label{lem:k'} The first column ${\cal K}' $ of
  \eqref{eq:Kdiag}
  is acyclic.
\end{lemma}
\pf We denote the maps by $\overline{\vi} $ and $\overline{\psi} $
respectively. To show that $\ker(\overline{\vi})= 0 $, consider an
element $u=\sum f_{kl} x_1^k x_2^l $ in the kernel. This is equivalent
to
\begin{equation} \label{eq:kervi} k f_{k,l} - m_1 f_{k+m_1, l+m_2} =
  0 = l f_{k,l} - m_2 f_{k+m_1, l+m_2} \text{\quad for $k \le -m_1 $
    and $l \le -m_2 $.}
\end{equation}
Now, $u \in L_1+L_2 $, so that $f_{kl}=0 $ if both $k $ and $l $ are
sufficiently negative and successive application of \eqref{eq:kervi}
gives $u \in K_{1,1} $.

\noindent Next, consider $(\overline{g}, \overline{h}) \!\in\!
\ker(\overline{\psi}) $, i.e. $\frac{\partial h}{\partial x_1}
- \frac{\partial g}{\partial x_2} + m_2 x_2^{-1} \alpha g - m_1
x_1^{-1} \alpha h \in K_{-2m_1, -2m_2} $. By definition, there are
numbers $N,M $ with $g,h \in K_{-M, -N} $. We want to solve
$\overline{\vi}(\overline{f}) = (\overline{g}, \overline{h}) $. To this
end, let $f_{kl}:=0 $ for all $k,l $ with $k<-M+m_1+1 $ and
$l<-N+m_2+1 $, but $(k,l) \neq (-M+m_1, -N+m_2) $. The desired
equation for $f $ induces the necessary equality of
$$
f_{-M+m_1, -N+m_2} = - \frac{1}{m_1} g_{-M-1, -N} = - \frac{1}{m_2}
h_{-M, -N-1} \ .
$$
Now, $\overline{\psi}(\overline{g}, \overline{h}) = 0 $ reads as
\begin{equation} \label{eq:coeff} (k+1) h_{k+1,l} - (l+1) g_{k, l+1}
  + m_2 g_{k+m_1, l+m_2+1} - m_1 h_{k+m_1+1, l+m_2} = 0
\end{equation}
for $k<-2m_1 $ and $l<-2m_2 $. We can choose representatives with
$g_{kl}=0 $ for $k \ge -m_1 $ or $l \ge -m_2+1 $, and $h_{kl}=0 $ for
$k \ge -m_1+1 $ or $l \ge -m_2 $. We can assume that $M>m_1+1 $ and
$N>m_2+1 $. Using \eqref{eq:coeff} for $k=-M-m_1 $, $l=-N-m_2 $ gives
the desired equality $m_2 g_{-M-1, -N} = m_1 h_{-M, -N-1} $ and
therefore the well-definedness of $f_{-M+m_1, -N+m_2} $.

In the same manner, we can solve the equation $k f_{kl} - m_1
f_{k+m_1, l+m_2} = g_{k-1,l} $ successively and get the coefficients
for the remaining indices in $A:= \{(k,l)|k \le 1, l \le 1 \} \setmin
\{(k,l)|k \le -M+m_1, l \le -N+m_2 \} $. We obtain 
$$
f_{kl}= - \frac{1}{m_1} g_{k-m_1-1, l-m_2} - \frac{k-m_1}{m_1^2}
g_{k-2m_1-1, l-m_2} - \ldots - \frac{(k-m_1) \cdots
  (k-rm_1)}{m_1^{r+1}} g_{k-rm_1, l-rm_2} 
$$
with $r:= \max \{ \frac{M}{m_1}, \frac{N}{m_2} \} $. From this it
follows that $f:=\sum_{(k,l) \in A} f_{kl} x_1^k x_2^l \in L_1+L_2 $
solves ${\rm pr_1}(\overline{\vi}(f)) = \overline{g} $ where ${\rm pr_1}
$ is the projection to the first direct summand. The same computations
give a solution $\widetilde{f} \in L_1+L_2 $ for the second direct
summand, i.e. ${\rm pr_2}(\overline{\vi}(\widetilde{f})) = \overline{h}
$. Using $\overline{\psi}(\overline{g}, \overline{h}) = 0 $ as above, we
easily see that the above calculations give the same coefficients
$f_{kl}= \widetilde{f}_{kl} $, so that we have found an $f \in
L_1+L_2 $ with $\widetilde{\vi}(f)=(\overline{g}, \overline{h}) $.

To see the surjectivity of $\overline{\psi} $, let $u \in
\mod{L_1+L_2}{K_{-2m_1, -2m_2}} $ be given, represented by $u \in
K_{-M,-N} $ for suitable $M,N $. Let $h:=0 $. We then have to find $g
$ whose coefficients satisfy
\begin{equation} \label{equ:coeff2} -l g_{kl} + m_2 g_{k+m_1, l+m_2}
  = u_{kl} \mbox{ for } k<-2m_1 \mbox{ and } l<-2m_2+1 \ .
\end{equation}
Especially, $lg_{kl} = m_2 g_{k+m_1, l+m_2} $ for $k< -M $ and $l<
-N-1 $, so that we can take $g_{kl}:=0 $ for $k<-M+m_1 $ and
$l<-N+m_2+1 $. Successively solving \eqref{equ:coeff2}, it follows
that 
$$
g_{kl}= \frac{1}{m_2} u_{k-m_1, l-m_2+1} + \frac{l-m_2}{m_2^2}
u_{k-2m_2, l-2m_2+1} + \ldots + \frac{(l-m_2) \cdots
  (l-am_2)}{m_2^{a+1}} u_{k-am_1, l-am_2+1}
$$
for $a \ge \max \{ \frac{M}{m_1}, \frac{N}{m_2} \} $. One obtains a
solution $g \in \mod{K_{-M+m_1, -N+m_2+1}}{K_{-m_1, -m_2+1}} $ and
hence the desired surjectivity of $\overline{\psi} $.
\\ \qed 

\medskip\noindent{\em 2. Step:} In order to calculate the cohomology
of $\mod{\Eess}{L_1+L_2} $, we thus have to do so for the middle
column ${\cal K} $ of \eqref{eq:Kdiag}. After transformation $x_i
\mapsto x_i^{-1} $ for $i=1,2 $, we have to consider the following
complex
\begin{equation} \label{equ:Hol} 0 \to \Hol \stackrel{D}{\lto}
  x_1^{m_1+1} x_2^{m_2} \Hol \oplus x_1^{m_1} x_2^{m_2+1} \Hol
  \stackrel{E}{\lto} x_1^{2m_1+1}x_2^{2m_2+1} \Hol \to 0 \ ,
\end{equation}
where $\Hol $ denotes the ring of power-series in two variables which
converge in the entire complex plane $\C^2 $. The maps are given as
follows. Define $P_{MN} := \{ \sum_{k,l \ge 0} f_{kl} x_1^k x_2^l \in
\Hol \mid f_{kl}=0 \mbox{ if } k>M \mbox{ and } l>N \} $.  Then
$x_1^{M+1} x_2^{N+1} \Hol \cong \mod{\Hol}{P_{MN}} $ and we put
$$
\begin{array}{l}
  D(u) \!\!:= \!\!-(x_1^2 \frac{\partial u}{\partial x_1} + m_1 x_1^{m_1+1}
  x_2^{m_2} u, x_2^2 \frac{\partial u}{\partial x_2} + m_2 x_1^{m_1}
  x_2^{m_2+1} u) {\rm mod} \, P_{m_1, m_2-1} \!\oplus\! P_{m_1-1, m_2}
  \\[0.2cm]
  E(\omega_1, \omega_2 ) := -x_1^2 \frac{\partial \omega_2}{\partial
    x_1} + x_2^2 \frac{\partial \omega_1}{\partial x_2} +m_2 x_2
  \alpha^{-1} \omega_1 - m_1 x_1 \alpha^{-1} \omega_2 \ {\rm mod} \
  P_{2m_1, 2m_2} \ .
\end{array}
$$
{\it Claim 1:} $\dim \ker D = (m_1, m_2) $.\\
Consider an element $u=\sum_{k,l \ge 0} u_{kl} x_1^k x_2^l \in \Hol $.
Then $u \in \ker D $ translates into the following condition on the
power series coefficients:
\begin{equation} \label{equ:kerD} k u_{kl} + m_1 u_{k-m_1, l-m_2} = 0
  \mbox{\qquad and \qquad} l u_{kl} + m_2 u_{k-m_1, l-m_2} = 0
\end{equation}
for all $k \ge m_1 $ and $l \ge m_2 $. It follows that for $u $ to be
in the kernel of $D $, it is necessary that its non-vanishing
coefficients $u_{kl} $ lie on the line ${\cal L} := \{ (k,l) \in \N_0
\times \N_0 \mid lm_1=km_2 \} $. Moreover, choosing values for the
coefficients $u_{kl} $ on this line in the region $0 \le k < m_1, 0
\le l < m_2 $ gives an element $u \in \ker D $ by means of
\eqref{equ:kerD} (the convergence of the solution so obtained is
easily seen). Claim 1 follows as the line ${\cal L} $ intersects the
integer lattice $\Z \times \Z $ in exactly $(m_1, m_2) $ points in
this region.

\noindent
{\it Claim 2:} $\dim {\rm coker} \ E=0 $.\\
Let $\eta \in x_1x_2 \alpha^{-2} \Hol $ be given. We have to find
$(\omega_1, \omega_2) \in x_1 \alpha^{-1} \Hol \oplus x_2 \alpha^{-1}
\Hol $ such that
\begin{equation} \label{equ:cokerE} -x_1^2 \frac{\partial
    \omega_2}{\partial x_1} + x_2^2 \frac{\partial \omega_1}{\partial
    x_2} + m_2 x_2 \alpha^{-1} \omega_1 - m_1 x_1 \alpha^{-1} \omega_2
  = \eta \ .
\end{equation}
We let $\omega_1 $ be arbitrary and write $\omega_2=
\exp(-x_1^{m_1}x_2^{m_2}) \cdot \rho $ with an element $\rho \in \Hol
$ yet to be determined. Then $(\omega_1, \omega_2) $ solves
\eqref{equ:cokerE} if and only if
$$
\frac{\partial \rho}{\partial x_1} = \exp(x_1^{m_1}x_2^{m_2}) \cdot
x_1^{-2} \left( x_2^2 \frac{\partial \omega_1}{\partial x_2} + m_2 x_2
  \alpha^{-1} \omega_1 - \eta \right) \ .
$$
The right hand side defines an element in $\Hol $ which can be
integrated in $x_1 $ direction, so that such a $\rho $ exists, proving
Claim 2.

\medskip\noindent{\em 3. Step:} We prove that the complex given in
\eqref{equ:Hol} has vanishing Euler characteristic. To this end,
consider the following diagram with exact rows:
\begin{equation} \label{eq:KKdiag}
  \begin{array}{rcccccccl}
    0 & \to & 0 & \to & \Hol & \stackrel{id}{\to} & \Hol & \to & 0 \\
    && \downarrow && \phantom{D} \downarrow D && \phantom{D} \downarrow D \\
    0 & \to & P_{m_1,m_2-1} \oplus P_{m_1-1, m_2} & \to & \Hol \oplus \Hol &
    \to & x_1 \alpha^{-1} \Hol \oplus x_2 \alpha^{-1} \Hol & \to & 0\\
    && \downarrow && \phantom{E} \downarrow E && \phantom{E} \downarrow E \\
    0 & \to & P_{2m_1, 2m_2} & \to & \Hol & \to & x_1x_2 \alpha^{-1} \Hol & \to & 0
  \end{array} \ .
\end{equation}
To complete this step, we prove the following
\begin{lemma} \label{lem:PFred} The operator
  $$
  E: P_{m_1, m_2-1} \oplus P_{m_1-1, m_2} \lto P_{2m_1, 2m_2}
  $$
  is Fredholm with index $-1 $. More precisely, one has $\dim \ker(E)
  = d+2 $ and $\dim {\rm coker}(E) = d+3 $, where $d:=(m_1, m_2) $
  denotes the greatest common divisor of $m_1 $ and $m_2 $.
\end{lemma}
\pf We first calculate the dimension of the cokernel. Let
$u:=\sum_{k,l} u_{kl} x_1^k x_2^l \in P_{2m_1,2m_2} $. We have to
solve $E(g,h)=u $. In terms of the coefficients in Laurent series
expansions this reads as
\begin{equation} \label{equ:coeff3} \psi_{kl} := -(k-1) h_{k-1,l} +
  (l-1) g_{k,l-1} + m_2 g_{k-m_1, l-m_2-1} - m_1 h_{k-m_1-1, l-m_2} =
  0
\end{equation}
for all $(k,l) \in A $, where $A:= \{ (k,l) \in \N_0 \times \N_0 \mid
k \le 2m_1 \mbox{ or } l \le 2m_2 \} $. Observe that in most cases
either the first two or the last two summands vanish, as
$g_{kl}=h_{kl}=0 $ for indices $k,l \not\in A $. There are several
different cases of pairs $(k,l) \in A $ to be considered:
\begin{enumerate}
\item $(k<m_1 $ or $l<m_2) $ and $(k>0 $ and $l>0) $: then
  \eqref{equ:coeff3} for the pair $(k,l) $ and the pair $(k+m_1,
  l+m_2) $ gives
  $$
  (k-1) h_{k-1,l} + (l-1) g_{k,l-1} = u_{kl} \mbox{ and }
  -m_1 h_{k-1,l} + m_2 g_{k,l-1} = u_{k+m_1, l+m_2} \ .
  $$
\item $k=0 $ and $l>0 $: Then \eqref{equ:coeff3} gives $\psi_{0,l}=
  (l-1)g_{0,l-1}=u_{0,l} $, $\psi_{m_1, l+m_2}=$
  $$
   = -(m_1-1) h_{m_1-1, m_2+l} + (m_2+l-1) g_{m_1,
     m_2+l-1} + m_2 g_{0,l-1} = u_{m_1, m_2+l}
  $$
  and $\psi_{2m_1, 2m_2+l} = m_2 g_{m_1, m_2+l-1} - m_1 h_{m_1-1,
    m_2+l} = u_{2m_1, 2m_2+l} $.
\item $k>0 $ and $l=0 $: in analogy to the case ii), this gives
  $\psi_{k,0}= -(k-1)h_{k-1,0}=u_{k,0} $, $\psi_{m_1+k, m_2}= $
  $$
  = -(m_1+k-1) h_{m_1+k-1, m_2} + (m_2-1) g_{m_1+k,
    m_2-1} - m_1 h_{k-1,0} = u_{m_1+k, m_2}
  $$
  and $\psi_{2m_1+k, 2m_2} = m_2 g_{m_1+k, m_2-1} - m_1 h_{m_1+k-1,
    m_2} = u_{2m_1+k, 2m_2+l} $.
\item $(k,l)=(0,0) $: Gives the equation $u_{0,0}=0 $.
\item $(k,l) \in \{ (m_1, m_2), (2m_1, 2m_2) \} $: Then
  \eqref{equ:coeff3} gives
  $$
\begin{array}{rl}
  \psi_{m_1, m_2} & = -(m_1-1) h_{m_1-1, m_2} + (m_2-1) g_{m_1,
    m_2-1} = u_{m_1, m_2} \mbox{ and } \\
  \psi_{2m_1, 2m_2} & = -m_1 h_{m_1-1, m_2} + m_2 g_{m_1, m_2-1}=
  u_{2m_1, 2m_2} \ .
\end{array}
$$
\end{enumerate}
By some simple matrix calculations one sees that case i) gives a $d
$-dimensional contribution to ${\rm coker}(E) $ if $m_1 \neq m_2 $ and
a $(d-1) $-dimensional contribution if $m_1=m_2 $; cases ii) and iii)
each give a one-dimensional contribution. Case iv) contributes with a
one-dimensional subspace, whereas case v) adds another dimension if
$m_1 \neq m_2 $ and no contribution for $m_1=m_2 $. Summing everything
up, gives the result
$$
\dim {\rm coker}(E)= d+3
$$
as claimed above. We remark that the convergence of the solutions we
obtain by the combinatorics of the Laurent coefficients above is an
easy exercise.

To compute the dimension of $\ker(E) $ we proceed in a similar manner,
looking at the cases i) -- v) above with $u=0 $. Again, one easily
sees that case i) give a $d $-dimensional subspace if $m_1 \neq m_2 $
and a $(d-1) $-dimensional one for $m_1=m_2 $, cases ii) and iii) each
add one dimension, case iv) gives no contribution to $\ker(E) $ and
case v) contributes with a one-dimensional subspace for $m_1=m_2 $ and
gives no solution for $m_1 \neq m_2 $. Hence, we have
$$
\dim \ker(E)= d+2 \ ,
$$
from which the lemma follows.
\\ \qed 

We are left with the task to calculate the index of the complex given
by the middle column of \eqref{eq:KKdiag}. To this end, again consider
the following diagram with exact rows:
\begin{equation} \label{eq:KKKdiag}
  \begin{array}{rcccccccl}
    0 & \to & \Hol & \stackrel{id}{\lto} & \Hol & \lto & 0 & \to & 0\\
    && \phantom{D'} \downarrow D' && \phantom{D} \downarrow D && \downarrow \\
    0 & \to & \Hol \oplus \Hol & \stackrel{-x_1^2, -x_2^2}{\lto} & \Hol \oplus \Hol
    & \lto & \mod{\Hol}{x_1^2 \Hol} \oplus \mod{\Hol}{x_2^2 \Hol} & \to & 0 \\
    && \phantom{E'} \downarrow E' && \phantom{E} \downarrow E &&
    \phantom{A} \phantom{A} \downarrow A \\
    0 & \to & \Hol & \stackrel{x_1^2x_2^2}{\lto} & \Hol & \lto &
    \mod{\Hol}{x_1^2x_2^2\Hol} & \to & 0 
  \end{array}
\end{equation}
where $A(\omega_1, \omega_2) := -x_1^2 \frac{\partial
  \omega_2}{\partial x_1} + x_2^2 \frac{\partial \omega_1}{\partial
  x_2} \, {\rm mod} \, x_1^1x_2^2 \Hol $. Obviously, $\ker(A) $ is the
$\C $-span of $(0,1), (0,x_2), (1,0), (x_1,0) $ and ${\rm coker}(A)=
{\rm span}_\C \{ 1, x_1, x_2, x_1x_2 \} $. It follows that $\dim \ker(A)= \dim
{\rm coker}(A)=4 $ and the Euler characteristic of the first and the
second column of \eqref{eq:KKKdiag} coincide. Thus vanishing of the
Euler-characteristic of \eqref{equ:Hol} follows from
\begin{lemma}
  The Euler-characteristic of $0 \!\to\! \Hol \stackrel{D'}{\to} \Hol
  \oplus \Hol \stackrel{E'}{\to} \! \Hol \to 0 $ equals $1 $.
\end{lemma}
\pf Decompose the operators as follows:
\begin{alignat}{1} \label{equ:stoer1}
  D'u & = \big( ( \frac{\partial}{\partial x_1},
    \frac{\partial}{\partial x_2}) +
    (m_1 x_1^{m_1-1} x_2^{m_2}, m_2 x_1^{m_1} x_2^{m_2-1}) \big) u
    \mbox{\qquad and} \\
  E'(\omega_1, \omega_2) & = (\frac{\partial \omega_2}{\partial x_1} -
  \frac{\partial \omega_1}{\partial x_2}) + (m_1 x_1^{m_1-1}x_2^{m_2}
  \omega_2 - m_2 x_1^{m_1} x_2^{m_2-1} \omega_1) \ .   \label{equ:stoer2}
\end{alignat}
For any $R>0 $, consider the space $B_R^{(k_1,k_2)} $ of all
holomorphic functions on $D_R \times D_R $ that are of type $C^{k_i}
$ with respect to the variable $x_i $ on $\overline{D}_R
\times \overline{D}_R $, where $D_R \subset \C $ denotes the open disc
around $0 $ with radius $R $ and $\overline{D}_R $ its closure. The
complex to be considered induces a complex $B^{(1,1)}_R
\stackrel{D'}{\lto} B^{(0,1)}_R \oplus B^{(1,0)}_R \stackrel{E'}{\lto}
B^{(0,0)}_R $.  The decompositions \eqref{equ:stoer1} and
\eqref{equ:stoer2} are decompositions of complexes in the sense that
taking either the first terms in \eqref{equ:stoer1} and
\eqref{equ:stoer2} or the second terms again give a complex. The
second terms in \eqref{equ:stoer1} and \eqref{equ:stoer2} consist of
compact operators and so the Euler-characteristic remains unchanged if
one omits these terms by Vasilescu's generalization of the well-known
compact perturbation theorem for Fredholm operators (cp. \cite{vas1}
and \cite{vas2}). The Euler-characteristic of the unperturbed complex
is easily seen to be $1 $, the only contribution coming from the
constant functions being the kernel of $(\frac{\partial}{\partial
  x_1}, \frac{\partial}{\partial x_2}) $ and the lemma follows by
taking the limit $R \to \infty $. Compare with \cite{mal4}, Theorem
1.4 for
the analogous arguments in the one-dimensional case.
\\ \qed  

\noindent
This lemma completes the proof of Theorem
\ref{thm:dRdim}. \\ \qed

\subsubsection{Non-degeneracy from the left}

Finally, we will prove non-degeneracy of the local pairing at a
crossing-point from the left. Note that perfectness from the left is
better accessible than perfectness from the right, the reason lying in the
difficulty of constructing sufficiently good 'test forms' $\omega \in
\Eess $ with $\nabla \omega \in L_1 + L_2 $, whereas 'test
cycles' $c \otimes \ve $ are easier to handle. The arguments used in
the proof are similar to the one-dimensional case of \cite{b-e}.

\begin{theorem} \label{thm:linksexakt} The pairing
  $$
  \Hdr^p(\mod{\Eess}{L_1 + L_2}) \times \Hrd_{p+2}(\Delta, \partial
  \Delta; D) \to \C \ , \quad p=0,1 \ ,
  $$
  is non-degenerate from the left, i.e. if $< [\omega], [c \otimes
  \ve] > = 0 $ for all $[c \otimes \ve] $, then $[\omega] = 0 $.
\end{theorem}
\pf We start with the case $p=0 $. Let $[\omega] \in
\Hdr^0(\mod{\Eess}{L_1+L_2}) $ be given. Let $\ve $ be a basis of
$\bve $ and $\vve $ denote the dual basis of $\bve^\vee $. Then
$\omega $ can be written as $\omega= a \cdot \ve $ with an analytic
function $a $ and
$$
\nabla \omega = \ve \otimes da = e \otimes
\eta_1 + e \otimes \eta_2 \in L_1^1+L_2^1 \ ,
$$
where $e $ denotes a meromorphic local basis of $L $. We have to
show that $a \ve \in L_1^0+L_2^0 $. Let $c $ be the radial
sheet $c=[0,p_1] \times [0,p_2] \subset
\Delta \cong D^2 \times D^2 $, with $p=(p_1, p_2) \in \partial \Delta
$.

First, consider the case that $\ve $ and $\vve $ are sections of
$L_1^0+L_2^0 $ in the notation used before. Since
$$
da = < \nabla \omega, \vve > = < e, \vve > \cdot \eta \ ,
$$
with $\eta=\eta_1+\eta_2 $, it follows that $da \in
\O_1^0+\O_2^0 $, hence also $a $ and therefore $a \ve \in
L_1^0+L_2^0 $, where as before $\O_\nu^p $ denotes the $p $-forms
meromorphic along $D_{\neq \nu} $ and arbitrary along $D_\nu $.

Next, assume that $\vve $ is rapidly decaying along $c $. Then by
assumption
\begin{equation} \label{equ:links} a(p) \vve(p) = (< \omega,
  \vve > \cdot \vve) (p) 
= \big( \Int_c < \nabla \eta_2, \vve
  > + \Int_{\partial c_1-D} < \eta_1, \vve > - \Int_{\partial
    c_2-D} < \eta_2, \vve > \big) \cdot \vve(p) \ ,
\end{equation}
with the decompositions $c=c_1+c_2 $ and $\nabla \omega= \eta_1 +
\eta_2 \in L_1^1+L_2^1 $ as in the definition of the pairing,
see Lemma \ref{lemma:alphadef}. Now,
$\nabla \eta_2 = - \nabla \eta_1 =  e \otimes \rho $ with a
meromorphic two-form $\rho $. With these notations, one has
$$
\Int_c < \nabla \eta_2, \vve > = \Int_c < e, \vve >
\rho \ .
$$
In local coordinates, the section $\vve $ is asymptotically equal to
$\exp(-kx_1^{-m_1}x_2^{-m_2}) $ times a meromorphic section of
$\bve^\vee $. Therefore, in order to understand the first term in
\eqref{equ:links} we have to study the behavior of
$$
e^{kp_1^{-m_1}p_2^{-m_2}} \cdot \int_{[0,p_1] \times [0,p_2]}
e^{-kx_1^{-m_1}x_2^{-m_2}} \cdot x_1^{-r_1} x_2^{-r_2} \, dx_1 \, dx_2
$$
for $(p_1, p_2) \to (0,0) $ and some $r_1, r_2 \in \Z $. Similar to
\cite{b-e} in one variable, substituting variables $y_i = x_i^{-1} $,
$q_i = p_i^{-1} $ and $u_i=y_i - q_i $, the latter integral reads as
$$
\Int_0^\infty \Int_0^\infty e^{k(-u_1^{m_1}u_2^{m_2}-f(u,q))} \cdot
(u_1+q_1)^{r_1-2} (u_2+q_2)^{r_2-2} \, du_1 \, du_2
$$
with a polynomial $f $ with positive coefficients. This integral has
at worst moderate growth as $(q_1,q_2) \to (\infty, \infty) $, so that
its contribution vanishes modulo $L_1+L_2 $.

In a similar manner, the second summand in \eqref{equ:links}
$$
\vve(p) \cdot \int_{\partial c_1-D} < \eta_1, \vve > = \big(
\int_{\{ p_1 \} \times [0, p_2]} < e, \vve > \eta_1
\big) \cdot \vve(p)
$$
leads us to study the integral
$$
e^{kp_1^{-m_1}p_2^{-m_2}} \cdot \Int_0^{p_2}
e^{-kp_1^{-m_1}p_2^{-m_2}} \cdot x_2^{-r_2} \cdot
\widetilde{\eta}(p_1) \, dx_2
$$
for $p_1 \to 0 $ with arbitrary $\widetilde{\eta} $. As before, this
has at most moderate growth for fixed $p_2 $ and $p_1 \to 0 $ and thus
lies in $L_2 $. The same argument shows that the third term in
\eqref{equ:links} vanishes modulo $L_1 $. It follows, that $a \ve
$ lies in $L_1^0+L_2^0 $ provided that $\vve $ is rapidly decaying
along $c $.

It remains to consider the case, where $\vve $ is rapidly increasing
along $c $, i.e. the radial sheet $c $ does not lie in any Stokes
bisector belonging to $\vve $. Then $\ve $ is rapidly
decaying along $c $. Again, by
$$
\nabla \omega = \ve \otimes da \in L_1^1+L_2^1
$$
it follows that $\ve \otimes da \in L_1^1+L_2^1 $. If we write
$\ve=\psi \cdot e$ with the analytic function $\psi $, this reads as
$\psi da \in \O_1^1+\O_2^1 $, i.e.
\begin{equation} \label{equ:psia} \psi \frac{\partial a}{\partial
    x_1} \in \O_1^0+\O_2^0 \mbox{\quad and \quad} \psi
  \frac{\partial a}{\partial x_2} \in \O_1^0+\O_2^0 \ .
\end{equation}
Now, $\psi $ has rapid decay along $c $ by assumption.  In order to
prove that \eqref{equ:psia} induces $\psi \cdot a \in \O_1^0+\O_2^0 $,
we have to show that for any two functions $g, a $ in the variables
$(x_1,x_2) $ such that $g $ has rapid decay and
$$
g \frac{\partial a}{\partial x_i} \in \O_1^0+\O_2^0 \mbox{ for } i=1,2
\ ,
$$
it follows that $ga \in \O_1^0+\O_2^0 $. We apply the mean value
theorem at a position $(x_1,x_2) $ after choosing a fixed point
$q=(q_1,q_2) $ with $0<|x_i|<|q_i| $ in order to find a point
$r=(r_1,r_2) $ in between such that
$$
g(x)a(x) = g(x) \cdot \big( a(q)+ (\frac{\partial a}{\partial
  x_1}(r),\frac{\partial a}{\partial x_2}(r)) \cdot \binom{h_1}{h_2}
\big) \ ,
$$
where $h=(h_1,h_2):=q-x $. By \eqref{equ:psia}, we can find functions
$\vi_\nu \in \O_\nu^0 $ such that $g(x)a(x) = g(x)a(q)+
\frac{g(x)}{g(r)} \cdot (\vi_1(r) + \vi_2(r)) $. Now, $g(x) $ is
rapidly decaying as $(x_1,x_2) \to (0,0) $ along $c $, so that there
is no growth contribution coming from the first term. As for the
second term, the function $\vi_1 $ has moderate growth in $x_2
$-direction and thus for fixed $x_1 $ and $r_1 $,
$$
| \frac{g(x)}{g(r)} \vi_1(r) | \le C_1 \cdot |r_2|^{-N} \le C_1 \cdot
|x_2|^{-N}
$$
proving moderate growth of the second term in $x_2 $-direction also.
The same argument applies for the third term, proving $ga \in
\O_1^0+\O_2^0 $.

Next, we consider the pairing $\Hdr^1(\mod{\Eess}{L_1+L_2}) \times
\Hrd_2(\Delta, \partial \Delta; D) \to \C $, i.e. the case $p=1 $. We
have the following Mayer-Vietoris sequence for the de Rham cohomology
of $\MM:= \mod{\Eess}{L_1+L_2} $ as well as the dual of the
corresponding sequence in rapid decay homology (see \eqref{equ:MV1})
$$\footnotesize\small
\begin{array}{ccccccccc}
  0 \ot & \Hdr^1(\pi^{-1}(0); \MM) & \ot & \bigoplus\limits_{j=1}^m
  \Hdr^0(T_{j,j+1}; \MM) & \ot & \bigoplus\limits_{j=1}^m \Hdr^0(T_j;
  \MM) & \ot & \Hdr^0(\pi^{-1}(0); \MM) & \ot 0 \\
  & \downarrow \alpha && \downarrow \beta && \downarrow \gamma && \downarrow\\
  0 \ot & \Hrd_3(\Delta, \partial)^\vee  & \ot &
  \bigoplus\limits_{j=1}^m \Hrd_2(Z_{j,j+1}, \partial)^\vee & \ot &
  \bigoplus\limits_{j=1}^m \Hrd_2(Z_j, \partial)^\vee & \ot &
  \Hrd_2(\Delta, \partial)^\vee & \ot 0
\end{array}\normalsize
$$
according to the decomposition of the torus $\pi^{-1}(0) =
\bigcup_{j=1}^m T_j $ with $T_j:=S^1 \times \mu_j $ where $S^1 $
decomposes into small enough intervals as $S^1=\bigcup_{j=1}^m \mu_j
$. The radial sheets with directions in $T_j $ are denoted by $Z_j $
as before, i.e. 
$$
Z_j := \bigcup\nolimits_{(\vt_1, \vt_2) \in T_j} \R^+_0
e^{i\vt_1} \times \R^+_0 e^{i\vt_2} \cap \Delta \ .
$$
Again, $T_{ij} :=
T_i \cap T_j $ and $Z_{ij} := Z_i \cap Z_j $. The vertical arrows are
induced from the local pairing. The arguments above show that $\beta $
and $\gamma $ are injective. Additionally, our knowledge about the
dimensions of the (co-)homology groups involved (Theorem
\ref{thm:rddim} and Theorem \ref{thm:dRdim}) then induce that $\gamma
$ is surjective, hence $\alpha $ is injective by the 5-lemma. This
proves perfectness from the
left for $p=1 $.
\\ \qed 

\subsection{Local pairing at a crossing-point involving one direction}
\label{sec:onedir}

To complete the investigation of the situation at a crossing point of
$D $, we have to consider the pairing involving the contributions from
the one-dimensional local strata of $D $. Let $\Delta $ denote a small
bi-disc around the crossing point $(0,0) \in D $, where we chose local
coordinates such that $D=\{ x_1x_2=0 \} = D_1 \cup D_2 $. We prove the
following

\begin{theorem} \label{thm:smooth} The local pairing
  \begin{equation} \label{equ:pairsmooth} \Hdr^{p-1}(\Delta;
    \mod{\Eess}{L_2}) \times \Hrd_p(\Delta \setmin D_2, \Delta
    \setmin D; D) \to \C
  \end{equation}
  for the local contribution of the connection along $D_1 $ at the
  crossing point is perfect.
\end{theorem}
\pf We omit $\Delta $ from the notation of de Rham
  cohomology.  Recall that $L $ has the elementary model $e^\alpha
  \otimes R $, where $\alpha(x)= x_1^{-m_\alpha} \cdot
  u(x) $ with $u(0) \neq 0 $.  

We now prove the analogue of Theorem \ref{thm:rddim}, namely
\begin{proposition} \label{prop:drdimj2} In the situation above, we
  have
  $$
  \dim \Hrd_\ast(\Delta \setmin D_2, \Delta \setmin D; D) = \left\{
    \begin{array}{ll}
      0 & \mbox{for } \ast \ge 3 \\
      m & \mbox{for } \ast=1,2 \ . \end{array} \right.
  $$
\end{proposition}
\pf We have the {\bf Stokes directions} in
  $\pi^{-1}(D_1^\times) \cong S^1 \times D_1^\times $:
$$
\Sigma_1 := \{ (\vt_1,x_2) \in \pi^{-1}(0,x_2) \mid -m_1 \vt_1 + {\rm
  arg}(u(0,x_2)) \in ( \frac{\pi}{2}, \frac{\pi}{2}) \} \ ,
$$
where $D_1^\times:= D_1 \setmin 0 $, and the fibration by {\bf Stokes
  sectors}, which reads in local coordinates as
$$
\S_1 := \bigcup_{(\vt_1,x_2) \in \Sigma_1} \R^+ e^{i \vt_1} \times
\{x_2 \} \ .
$$
Note, that one may regard $\S_1 $ as a subset of the disc bundle
associated to the normal bundle $DN(D_1^\times) $ of $D_1^\times $ in
$X $, which we identify with $\Delta \setmin D_2 $.  In each fiber
over $(0,x_2) \in D_1^\times $ it is determined by the Stokes sectors
of $e^\alpha $, i.e. the sectors, where $e^\alpha $ has rapid decay
for $x_1 \to 0 $.

Decomposing $S^1 $ into small enough intervals and thereby $\Delta
\setmin D_2$ into the corresponding fibration by sectors, we can apply
Mayer-Vietoris sequences in the same way as we did in the proof of
Proposition \ref{prop:rdcrossiso}. With the same arguments as before,
the statement analogous to Proposition \ref{prop:rdcrossiso} follows,
namely a decomposition as ordinary singular homology groups:
$$
\Hrd_\ast(\Delta \setmin D_2, \Delta \setmin D;D) \cong H_\ast(V_1
\cap \S_1, \partial V_1 \cup D) \ ,
$$
where $\S_1 $ denotes the fibration by Stokes sectors
corresponding to $e^\alpha $ and $V_1 $ denotes a small tubular
neighborhood of $D_1^\times := D_1 \setmin 0 $.

It remains to show, that
\begin{equation} \label{equ:stokes1} \dim H_\ast(V_1 \cap
  \Sigma_1, \partial V_1 \cup D; D) = \left\{ \begin{array}{ll}
      0 & \mbox{for } \ast \ge 3 \\ m_\alpha & \mbox{for } \ast=1,2
    \end{array} \right.
\end{equation}

The topological situation looks as follows. Consider the fiber
$DN(D_1^\times)_x $ of the normal disc bundle of $D_1^\times $ at some
point $x:=(0,x_2) \in D_1^\times $.  We denote by $\sigma^\alpha
\subset DN(D_1^\times)_x $ the Stokes sectors inside this disc, i.e.
the intersection $\sigma^\alpha = \Sigma^\alpha_1 \cap
DN(D_1^\times)_x $. Now, we can retract the tubular neighborhood $V_1
$ of $D_1^\times $ to the trivial bundle with base space being a
circle $S^1 $ around $D_2 $ and with fiber $DN(D_1^\times)_x $, as
well as the Stokes sectors $\sigma_\alpha $ in this fiber to the union
of $m_\alpha $ radii, one for each Stokes sector. Now, we have to
identify all points in the boundary $\partial V_1 $ as well as those
in $D $. In the fiber $DN(D_1^\times)_x $, this result in a wedge of
$m $ circles, and we finally are left with a trivial bundle
over the circle $S^1 $ with fiber the wedge $\bigvee_{m_\alpha} S^1 $,
where we still have to identify the 'zero'-section $S^1 \times \{ pt
\} \subset D_1^\times $, $\{ pt \} $ being the base point of the wedge
of the $m_\alpha $ circles. Thus,
$$
H_\ast(V_1^\times \cap \Sigma^\alpha_1, \partial V_1 \cup D) \cong
H_\ast \big( \mod{(S^1 \times \bigvee\limits_{m_\alpha} S^1)}{(S^1
  \times \{ pt \}}) \big) \ .
$$
Using the K\"unneth isomorphism and the effect of collapsing the
one-cell $S^1 \times \{ pt \} $, one easily sees that the latter
singular homology spaces obviously have the desired dimensions as in
\eqref{equ:stokes1}.
\\ \qed 

Next, we calculate the dimension of the corresponding de Rham
cohomology groups from the pairing \eqref{equ:pairsmooth}:
\begin{proposition} \label{prop:rddimj2} In the situation from above,
  we have
  $$
  \dim \Hdr^\ast(\mod{\Eess}{L_2}) = \left\{ \begin{array}{ll} 0 &
      \mbox{for } \ast \neq 0,1 \\
      m & \mbox{for } \ast = 0,1 \ . \end{array} \right.
  $$
\end{proposition}
\pf As before, the regular singular part $R $ gives no
  contribution to the dimension, so that we can assume it is trivial.
  Again, we write $\alpha=x_1^{-m} \cdot u_\alpha(x) $ with
  $u_\alpha(0) \neq 0 $.  After a local transformation by
  $u_\alpha^{-1} $, we can assume $u_\alpha \equiv 1 $. We thus have
  to consider the complex
$$
\mod{\Eess}{L_2} \stackrel{\vi}{\lto} \mod{\Eess}{L_2} \oplus
\mod{\Eess}{L_2} \stackrel{\psi}{\lto} \mod{\Eess}{L_2}
$$
with
\begin{equation} \label{equ:vipsi} \vi(u):= (\frac{\partial
    u}{\partial x_1} - m x_1^{-m-1} u, \frac{\partial u}{\partial
    x_2}) \, \mbox{and} \, \psi(\omega_1, \omega_2) := \frac{\partial
    \omega_2}{\partial x_1} - m x_1^{-m-1} \omega_2 - \frac{\partial
    \omega_1}{\partial \omega_2} \ .
\end{equation}

Now, let $\o $ denote the ring of holomorphic germs in one variable $z
$ at $0 $ and $\M $ denote the meromorphic germs. Let
$\rho:\mod{\o}{\M} \to \mod{\o}{\M} $ be the map induced form the
one-variable connection $\nabla 1:= z^{-m-1} dz $, i.e.
$\rho(f)=f'(z)-m z^{-m-1} f(z) $. Consider the diagram
$$
\begin{array}{ccccccc} f & \in & \mod{\o}{\M} & \stackrel{\rho}{\lto}
  & \mod{\o}{\M} & \ni & g \\
  \downarrow && \downarrow \beta && \downarrow \gamma && \downarrow\\
  f(x_1) & \in & \mod{\Eess}{L_2} & \stackrel{\vi}{\lto} &
  \mod{\Eess}{L_2} \oplus \mod{\Eess}{L_2} & \ni & (g(x_1),0) \ .
\end{array}
$$
We claim, that $\beta $ induces an isomorphism of the kernels of the
horizontal arrows. To this end, consider $u=\sum_{i,j} u_{ij}x_1^i
x_2^j $. Then $u \, {\rm mod} L_2 \in \ker \vi $ if and only if
$$
{\rm i)} \quad \frac{\partial u}{\partial x_2}=: \eta \in L_2
\mbox{\quad and \quad} {\rm ii)} \quad \frac{\partial u}{\partial x_1}
-m x_1^{-m-1} u \in L_2
$$
are both meromorphic in $x_1 $-direction. Assume that $\eta = \sum_{i
  \ge -N} \sum_j \eta_{ij} x_1^i x_2^j $, then by i) it follows that
$$
u_{i,j+1} = \frac{1}{j+1} \cdot \eta_{i,j} \mbox{\quad for all } j
\neq -1 \ ,
$$
especially $u_{i,j}=0 $ for all $j \neq 0 $ and all $i<-N $. Let
$\widetilde{u}:= \sum_{i \ge -N} \sum_j u_{ij}x_1^i x_2^j \in L_2 $,
then
$$
u(x_1,x_2) - \widetilde{u}(x_1,x_2) = \sum_{i<-N} u_{i,0} x_1^i =:
w(x_1) \in {\rm im}(\beta) \ .
$$
Now, $[u]=[\beta(w)] \in \mod{\Eess}{L_2} $ and obviously $[u] \in
\ker \vi \Leftrightarrow [w] \in \ker \rho $. From the computation of
the one-variable case (\cite{mal4}), we deduce $\dim \ker \vi = \dim
\ker \rho = m $.

Next, we want to compute $H^2(\mod{\Eess}{L_2})=0 $, i.e. for given
$\eta \in \Eess $, we have to find $(\omega_1, \omega_2) \in \Eess
\oplus \Eess $ and $\nu \in L_2 $, such that
$$
\frac{\partial \omega_2}{\partial x_1} - \frac{\partial
  \omega_1}{\partial x_2} - m x_1^{-m-1} \omega_2 = \eta+\nu \ .
$$
We can do so by setting $\omega_1:=0 $. Writing
$\omega_2=e^{-x_1^{-m}} \cdot f $, we obtain a solution if and only if
$\frac{\partial f}{\partial x_1} = e^{x_1^{-m}}(\eta+\nu) $, which is
a question of vanishing of the residue in $x_1 $-direction (i.e. the
residue of the function in the variable $x_1 $ for fixed $x_2 $). A
simple calculation shows that this can be achieved by
$$
\nu:= - \sum_j \eta_{-mk-1,j} \cdot x_2^j x_1^{-mk-1} \in L_2 \ .
$$
It remains to prove that $\dim \Hdr^1(\mod{\Eess}{L_2})=m $. Consider
an element $[(\omega_1, \omega_2)] \in \mod{\ker(\psi)}{{\rm im}(\vi)}
$.  We can write $\omega_2= \sum_j a_j(x_1) \cdot x_2^j $ and then
find some $u \in \Eess $ such that $\frac{\partial u}{\partial x_2}=
\sum_{j \neq -1} a_j(x_1)x_2^j $. Modulo ${\rm im}(\vi) $ we may
therefore assume that $\omega_2=a(x_1) \cdot x_2^{-1} $. We write
$\omega_1 $ in the form $\omega_1= \sum_j b_j(x_1)x_2^j $. Since
$[\omega_1, \omega_2] \in \ker(\psi) $, there exists an $\eta \in L_2
$, written as $\eta= \sum_j \eta_j(x_1) \cdot x_2^j $, such that
$$
\frac{\partial \omega_2}{\partial x_1} - mx_1^{-m-1} \omega_2 = \eta +
\frac{\partial \omega_1}{\partial x_2} \ ,
$$
which reads as $(a'(x_1) -m x_1^{-m-1} a(x_1)) x_2^{-1} = \sum_j
(\eta_j(x_1)+ (j+1)b_j(x_1)) x_2^j $.  It follows that $b_j(x_1)=
\frac{-1}{j+1} \eta_j(x_1) \in \M $ for any $j \neq -1 $, so that
modulo $L_2 $, we may again assume that $\omega_1=b(x_1) x_2^{-1} $
with $b(x_1):=b_{-1}(x_1) $. Furthermore the function $a(x_1) $ must
satisfy
$$
a'(x_1) -mx_1^{-m-1} a(x_1) = \eta_{-1}(x_1) \in \M \ .
$$
It follows that $\frac{\partial \omega_2}{\partial x_1} -m x_1^{-m-1}
\omega_2 \in L_2 $ and thus
$$
-b(x_1)x_2^{-2} = \frac{\partial \omega_1}{\partial x_2} =
\frac{\partial \omega_2}{\partial \omega_1} -m x_1^{-m-1} \omega_2 -
\eta \in L_2 \ .
$$
But then $\omega_1 \in L_2 $ and we have obtained
$$
\Hdr^1(\mod{\Eess}{L_2}) = \{ [(0,a(x_1)x_2^{-1})] \mid
a'(x_1)-mx_1^{-m-1} a(x_1) \in \M \} \ .
$$
From the theory in the case of one variable, we know that its
dimension is $m $ as we wanted to prove.  \\ \qed 

Thus, we have proved that the (co-)homology spaces in
\eqref{equ:pairsmooth} have the same dimension. It remains to prove
perfectness from one side:
\begin{proposition} \label{prop:smoothlinks} The local pairing
  \eqref{equ:pairsmooth} is non-degenerate from the left, i.e. assuming $<
  [\omega], [c \otimes \ve] > = 0 $ for all $[c \otimes \ve] $, then
  $[\omega]=0 $.
\end{proposition}
\pf The proof uses literally the same arguments as the one for the
corresponding statement for the pairing at a crossing-point including
both directions $D_1 $ and $D_2 $ (Theorem \ref{thm:linksexakt}) and
is therefore omitted here.
\\ \qed 

\noindent
The proposition completes the proof of Theorem \ref{thm:smooth}.
\\ \qed

\subsection{Local pairing at a smooth point}
\label{sec:smoothcomp}

It remains to study the local pairing at a smooth point of $D $, i.e.
we may now assume that in local coordinates $D:=D_1=\{ x_1=0 \} $ and
we have to consider a small bi-disc $\Delta $ around the smooth point
$(0,0) \in D $. Our task is to prove perfectness of the local pairing
\begin{equation} \label{equ:wirklsmooth} \Hdr^p(\Delta, \Delta
  \setmin D_1;D_1) \times \Hrd_{p+1}(\Delta; \mod{\Eess}{\Emero}) \lto
  \C \ .
\end{equation}

\begin{theorem} \label{thm:wirklsmooth} The pairing
  \eqref{equ:wirklsmooth} is perfect.
\end{theorem}
\pf The proves of the following propositions is very similar to the
proves of Proposition \ref{prop:drdimj2}, \ref{prop:rddimj2} and
\ref{prop:smoothlinks}, we will only briefly mention on the necessary
changes. Let $e^\alpha \otimes R $ be the elementary model with
$\alpha(x)=x_1^{-m_\alpha} u(x) $ as before.  
\begin{proposition} One has 
  $
  \dim \Hrd_\ast(\Delta, \Delta \setmin D_1; D_1) = \left\{
    \begin{array}{ll} 0 & \mbox{for } \ast \ge 2 \\ m_\alpha &
      \mbox{for } \ast=1 \end{array} \right. \ .
  $
\end{proposition}
\pf As in Proposition \ref{prop:rddimj2}, the rapid decay homology
groups decomposes as the ordinary singular homology groups of the
Stokes regions in the tubular neighborhood of $D_1 $, i.e. the normal
disc bundle of $D_1 $ modulo boundary and basis $D_1 $. Here, the
topological situation is even simpler as before, as we can retract the
base of the bundle, namely $D_1 $, to a point and end up with the
Stokes radii in the fiber over a chosen point $(0,x_2) \in D_1 $.
Modulo boundary and modulo base point $x $, this is just the wedge of
$m_\alpha $ circles. Thus
$$
\dim \Hrd_\ast(\Delta, \Delta \setmin D_1; D_1) = \dim H_\ast(
\bigvee_{m_\alpha} S^1) \ .
$$
\\ \qed 
\begin{proposition} One has $
  \dim \Hdr^\ast(\Delta; \mod{\Eess}{\Emero}) = \left\{
    \begin{array}{ll} 0 & \mbox{for } \ast \neq 0 \\ m_\alpha &
      \mbox{for } \ast=0
    \end{array} \right. \ .
  $
\end{proposition}
\pf The proof for degree $0 $ and $2 $ is nearly literally the same as
in Proposition \ref{prop:drdimj2}: We have to consider
$$
\mod{\Eess}{\Emero} \stackrel{\vi}{\lto} \mod{\Eess}{\Emero} \oplus
\mod{\Eess}{\Emero} \stackrel{\psi}{\lto} \mod{\Eess}{\Emero}
$$
with $\vi $ and $\psi $ defined as in \eqref{equ:vipsi}. Now, we
proceed in the same way considering the Laurent expansions. In the
case here, these expansions have no polar part in $x_2 $-direction,
but the arguments used for degree $0 $ and $2 $ still remain valid and
give the desired result.

The non-existence of polar parts in $x_2 $-direction, however, induces
vanishing of $\Hdr^1(\mod{\Eess}{\Emero}) $: Let $[(\omega_1, \omega_2)]
\in \ker(\psi) $. Now, $\omega_2 $ has no polar part in $x_2
$-direction, hence no residue in $x_2 $ (with fixed $x_1 $), so that
$$
\frac{\partial u}{\partial x_2} = \omega_2
$$
has a solution in $\Eess $. Therefore, modulo the image of $\vi $, we
can assume that $\omega_2=0 $. But then $\omega_1 $ has to fulfill the
equation
$$
\frac{\partial \omega_1}{\partial x_2} = \eta = \sum_{j \ge 0}
\eta_j(x_1) x_2^j
$$
for some $\eta \in \Emero $, i.e. with meromorphic functions $\eta_j
$. But this forces $\omega_1 $ to be meromorphic also, hence
$[(\omega_1, \omega_2)]= 0 \in \Hdr^1(\mod{\Eess}{\Emero}) $.
\\ \qed 

\begin{proposition}
  The pairing \eqref{equ:wirklsmooth} is non-degenerate from the left.
\end{proposition}
\pf The proof uses the same arguments as in Theorem
\ref{thm:linksexakt} and is therefore omitted here.
\\ \qed 

\noindent
We thus have proved local perfectness in all three situations to be
considered (cp. Definition \ref{def:locperf}). This finally completes
the proof of the main result (Theorem \ref{thm:main}).

\bigskip\noindent
{\small {\bf Acknowledgements:} I am grateful to Spencer Bloch
  and H\'el\`ene Esnault for introducing me to the subject of
  irregular singular connections, their permanent interest in my work
  and several fruitful discussions. Large parts of the work presented
  here were achieved while I stayed at the University of
  Chicago supported by the Feodor-Lynen program of the Alexander von
  Humboldt-Stiftung. I thank the Humboldt-Stiftung for offering me
  this opportunity as well as the University of Chicago for its
  hospitality and especially Spencer Bloch for the many discussions on
  this work during my stay. I also thank Stefan Bechtluft-Sachs,
  Uwe Jannsen, Byungheup Jun and Tomohide Terasoma for various helpful
  conversations. Finally, I thank C. Sabbah for numerous remarks on a
  previous version of this work and I am especially indebted to him
  for pointing out to me the importance of {\em good} formal models
  for connections on surfaces (see Definition \ref{def:good}).}


\begin{thebibliography}{3}
\bibitem{balda} F. Baldassarri, {\it Comparaison entre la
    cohomologie alg\'ebrique et la cohomologie $p $-adic rigide \`a
    coefficients dans un module diff\'erentiel I. (Cas de courbes)},
  Invent. Math. 87 (1987), 83 -- 99

\bibitem{bbe} A. Beilinson, S. Bloch, H. Esnault, {\it $\ve
    $-factors for Gauss-Manin determinants}, Mosc. Math. J. 2(2002),
  no. 3, 477 -- 532
  
\bibitem{b-e} S. Bloch, H. Esnault, {\it Homology for irregular
    connections}, Journal de Th\'eorie des Nombres de Bordeaux 16
  (2004), 65 -- 78
  
\bibitem{b-e:irreg} S. Bloch, H. Esnault, {\it Gauss-Manin
    determinant connections and periods for irregular connections},
  GAFA 2000 (Tel Aviv 1999), Geom. Funct. Anal. 2000, Special Volume,
  Part I, 1 -- 31
  
\bibitem{b-e:gm} S. Bloch, H. Esnault, {\it Gauss-Manin
    determinants for rank $1 $ irregular connections on curves. With
    an appendix in French by P. Deligne}, Math. Ann. 321(2001), no. 1,
  15 --87
  
\bibitem{borel} A. Borel, {\it Algebraic ${\cal D} $-modules},
  Perspectives in Math. 2, Academic Press, Boston 1987
  
\bibitem{deligne} P. Deligne, {\it Equations diff\'erentielles \`a
    points singuliers r\'eguliers}, LNS 163, Springer-Verlag, Berlin
  Heidelberg, 1970
  
\bibitem{groth} A. Grothendieck, {\it On the de Rham cohomology of
    algebraic varieties}, Publ. Math. IHES 29 (1966), 93 -- 103
  
\bibitem{haraoka} Y. Haraoka, {\it Confluence of cycles for
    hypergeometric functions on $Z_{2,n+1} $}, Trans. AMS 349, no. 2
  (1997), 675 --712
  
\bibitem{katz} N. Katz, {\it On the calculations of some
    differential Galois groups}, Inv. Math. 87 (1987), 13 -- 61

\bibitem{kihata} H. Kimura, Y. Haraoka, K. Takano, {\it The
    generalized confluent hypergeometric functions}, Proc. Japan
  Acad. Ser A 68 (1992), 290 -- 295
  
\bibitem{majima} H. Majima, {\it Asymptotic analysis for
    integrable connections with irregular singular points}, LNS 1075,
  Springer-Verlag, Berlin, Heidelberg, 1984
  
\bibitem{mal1} B. Malgrange, {\it Equations diff\'erentielles
    \`a coefficients polynomiaux}, Prog. in Math. 96, Birkh{\"a}user,
  Basel, Boston, 1991
  
\bibitem{mal2} B. Malgrange, {\it La classification des
    connexions irr\'eguli\`eres \`a une variable}, in 'Seminaire
  E.N.S. Math\'ematique et Physique', Progr. in Math. 37,
  Birk{\"a}user, Basel, Boston 1983, 381 -- 399
  
\bibitem{mal3} B. Malgrange, {\it Ideals of differentiable
    functions}, Oxford University Press, 1966

\bibitem{mal4} B. Malgrange, {\it Sur les points singuliers des
    equations diff\'erentielles}, L'enseignment math., t. XX (1974),
  fasc. 1-2, 147 -- 176
  
\bibitem{mebkhout} Z. Mebkhout, {\it Le th\'eor\`eme de
    comnparaison entre cohomologies de de Rham d'une vari\'et\'e
    alg\'ebrique complexe et le th\'eor\`eme d'existence de Riemann},
  Publ.Math. IHES 69(1989), 47-89
  
\bibitem{sabbah1} C. Sabbah, {\it Equations diff\'erentielles \`a
    points singuliers irr\'eguliers et ph\'enom\`ene de Stokes en
    dimension 2}, Ast\'erisque 263, Soc. Math. de France, 2000
  
\bibitem{sabbah2} C. Sabbah, {\it Equations diff\'erentielles \`a
    points singuliers irr\'eguliers en dimension 2}, Ann. Inst.
  Fourier (Grenoble) 43 (1993), 1619 -- 1688
  
  
\bibitem{terasoma} T. Terasoma, {\it Confluent hypergeometric
    functions and wild ramification}, J. of Algebra 185 (1996), 1 --
  18
  
\bibitem{vas1} F.-H. Vasilescu, {\it Stability of the index of a
    complex of Banach spaces}, J. Operator Theory 2(1979), 247 -- 275
  
\bibitem{vas2} F.-H. Vasilescu, {\it The stability of the Euler
    characteristic for Hilbert complexes}, Math. Ann. 248(1980), 109
  -- 116

\end{thebibliography}
\end{document}